\documentclass[a4paper]{article}
\pdfoutput=1

\usepackage{IEK10} 
\usepackage{natbib}

\makeatother
\usepackage{amssymb}
\usepackage{xcolor}
\usepackage{mathtools, cuted}
\usepackage{tabularx}
\usepackage{eurosym}
\usepackage[utf8]{inputenc}

\renewcommand{\min}[1][]{
	\ifthenelse{\isempty{#1}}{\operatorname{min}}{\ensuremath{\underset{#1}{\text{min}\,}}}
}

\def\IEK10{
Institute of Energy and Climate Research -- Energy Systems Engineering (IEK-10),
  Forschungszentrum Jülich GmbH,
  Wilhelm-Johnen-Straße,
  52425 Jülich,
  Germany
}
\def\LTT{
  Institute of Technical Thermodynamics, 
  RWTH Aachen University,
  Schinkelstraße 8,
  52062 Aachen,
  Germany
}
\def\ETH{
    Energy \& Process Systems Engineering,
    Department of Mechanical and Process Engineering, 
    ETH Zürich,
  Tannenstrasse~3,
  8092 Zürich,
  Switzerland
}
\def\CUR{
    Energy Trading and Dispatching, Currenta GmbH \& Co. OHG,
    Kaiser-Wilhelm-Allee 80,
51373 Leverkusen,
  Germany
}

\newcommand{\mytitle}{Where to Market Flexibility? Optimal Participation of Industrial Energy Systems in Balancing-Power, Day-Ahead, and Continuous Intraday Electricity Markets}

\newcommand{\affil}{
  \begin{itemize}[leftmargin=3mm, itemsep=0mm]
    \item[$^a$]\ETH
    \item[$^b$]\LTT
    \item[$^c$]\CUR
    \item[$^d$]\IEK10
  \end{itemize}
}

\def\firstAuthor{Niklas Nolzen}
\newcommand{\myauthor}{\firstAuthor$^{a,b}$, Alissa Ganter$^{a,b}$, Nils Baumgärtner$^{c}$, Ludger Leenders$^{a}$, André Bardow$^{a,d,*}$}
\newcommand{\correspondingauthor}{A. Bardow, \ETH \newline E-mail: abardow@ethz.ch}
\author{\myauthor}

\usepackage[
  colorlinks,
  linkcolor=blue,
  citecolor=blue,
  urlcolor=blue,
  pdftitle={\mytitle},
  pdfauthor={\firstAuthor}
]{hyperref}
\usepackage[capitalise, nameinlink]{cleveref}
\crefname{table}{Tab.}{Tab.}

\newcommand{\setpgfexternalcounter}[1]{
  \makeatletter%
  \pgfkeysgetvalue{/tikz/external/figure name}\myexternalname
  \expandafter\gdef\csname c@tikzext@no@\myexternalname\endcsname{#1}%
  \makeatother
}

\begin{document}

  \thispagestyle{firststyle}

  \begin{center}
    \begin{large}
      \textbf{\mytitle}
    \end{large} \\
    \myauthor
  \end{center}

  \vspace{0.5cm}

  \begin{footnotesize}
    \affil
  \end{footnotesize}

  \vspace{0.5cm}

  \begin{abstract}
The rising share of volatile renewable generation increases the demand for flexibility in the electricity grid.  
Flexible capacity can be offered by industrial energy systems through participation on either the continuous intraday, day-ahead, or balancing-power markets.
Thus, industrial energy systems face the problem of where to market their flexible capacity.
Here, we propose a method to integrate trading on the continuous intraday market into a multi-market optimization for flexible industrial energy systems.
To estimate the intraday market revenues, we employ option-price theory.
Subsequently, a multi-stage stochastic optimization determines an optimized bidding strategy and allocates the flexible capacity.
The method is applied to a case study of a multi-energy system showing that coordinated bidding in all three considered markets reduces cost most.
A sensitivity analysis for the intraday market volatility reveals changing market preferences, thus emphasizing the need for multi-market optimization. 
The proposed method provides a practical decision-support tool in short-term electricity and balancing-power markets.
\end{abstract}

\vspace{0.5cm}


\vspace{0.75cm}

\begin{flushleft}
  \leavevmode{\parindent=15mm\indent}
  \textbf{Highlights:}
  \begin{itemize}[leftmargin=20mm]
    \item Method optimizes a coordinated bidding strategy for all markets
    \item Integration of the continuous intraday market into a multi-market optimization
    \item Option-price theory estimates profit from trading in the continuous intraday market
    \item Multi-market participation leads to lowest costs for multi-energy systems 
  \end{itemize}
\end{flushleft}
\vspace*{5mm}

  \newpage

\section{Introduction}
\label{sec:Introduction}
\subsection{The challenges and opportunities of multiple markets}
The recent rise of variable renewable generation increases supply uncertainty in the electricity grid.
Thus, more short-term system balancing is needed to ensure grid stability.
The growing need for grid balancing is addressed by various measures, such as international grid control cooperation for balancing power \citep{Ocker.2017} and intraday markets \citep{Koch.2019}.

Continuous intraday markets have been introduced all over Europe in the last decade \citep{Ocker.2020}.
The intraday market continuously trades electricity until gate-closure, the market's closing time, usually shortly before delivery.
The continuous trading scheme allows adjusting for unexpected deviations in electricity supply and demand as more precise weather forecasts are available shortly before delivery.
Thereby, the continuous intraday market has become an important measure to reduce the deviation between supply and demand before gate-closure \citep{Koch.2019}.

After gate-closure of the continuous intraday market, remaining deviations are compensated with balancing power:
If the overall supply is lower than demand, the transmission system operator requests positive balancing power, and vice versa for negative balancing power \citep{Ocker.2017}.

The ba\-lan\-cing-power market and the continuous intraday market offer different marketing opportunities for flexible capacity:
On the continuous intraday market, flexible capacity can exploit price volatilities by asset-backed trading \citep{Lohndorf.2021}.
Asset-backed trading continuously trades the flexible capacity based on the real-time electricity price and the marginal electricity costs:
If electricity prices rise above marginal electricity costs, electricity is sold to the market.
If electricity prices fall below marginal electricity costs, electricity is purchased from the market.
Eventually, the sum of all purchases and sales determines the amount of electricity provided by the flexible capacity.

On the balancing-power market, flexible capacity commits positive and/or negative balancing power for a certain time period in the future, e.g., the next day or the following week.
Depending on the type of balancing power and the market design, the remuneration consists of the capacity price for the provision of balancing-power capacity and/or the energy price for the request of balancing power \citep{Barbero.2020}.

Thus, flexible capacity can be monetized on either the balancing-power market or the continuous intraday market.
However, bidding on short-term electricity and balancing-power markets leads to a sequential decision-making process.
E.g., in Germany, participation in the balancing-power market is followed by a tender in the day-ahead market.
Finally, after clearing of the day-ahead market, continuous trading starts and continues until shortly before delivery.
During delivery, the transmission system operator may request balancing power.

This sequential decision-making process couples the decisions:
The availability of electricity and flexibility to be traded on one market depends on the commitments to the other markets.
Hence, if flexible capacity is monetized on one market, this capacity is no longer available to other markets.

Optimal participation in sequential electricity markets thus requires coordinated bidding \citep{Lohndorf.2022}.
For coordinated bidding, the value of the opportunity to trade on the continuous intraday market needs to be considered while deciding on participation in the balancing-power market and day-ahead market \citep{Klaboe.16.06.201320.06.2013}.
As a result, flexibility providers face the question of where to best market their flexibility.

\subsection{Literature review on multi-market participation}
Coordinated bidding considers participation in multiple markets in one single optimization problem \citep{Aasgard.2019}.
Participation in the balancing-power market and day-ahead market is usually modeled as multi-stage stochastic optimization to account for uncertainties.
Since the early work of \citet{Swider.2007} and \citet{Boomsma.2014}, models have incorporated the market rules in more detail, particularly regarding the balancing-power market.

\citet{Muche.2016} consider average capacity and energy prices, thereby neglecting optimized bidding decisions on the balancing-power market.
\citet{Kumbartzky.2017} and \citet{Schafer.2019} optimize only the capacity price bids, whereas \citet{Leenders.2020} model the request of balancing power, thus only optimizing the energy price bids.
\citet{Bohlayer.2020b} model the market participation in the balancing-power and day-ahead market in a multi-stage stochastic optimization model that considers both the acceptance and the request of balancing power.
The model includes the most important aspects of the balancing-power market but disregards the continuous intraday market. 
Hence, the reviewed publications \citep{Swider.2007} -- \citep{Bohlayer.2020b} do not consider the trade-offs in allocating flexibility on all short-term markets.

Only few studies also consider intraday markets as counterpart to the balancing-power market.
\citet{Dowling.2017} model the participation in intraday, day-ahead, and balancing-power markets.
The approach assumes simultaneous trading on all markets and, thus, does not consider the sequential decision-making process required, e.g., for European markets.	
\citet{Kraft.2022} model a power plant portfolio that participates in the German balancing-power market, day-ahead market, and continuous intraday market.
This approach accounts for price uncertainties and price correlations among the markets.
However, the request and the remuneration of balancing power are not explicitly modeled, thus neglecting a potential revenue stream of the balancing-power market.
Additionally, the continuous intraday market is modeled as one stage with a single recourse decision.
Hence, both \citet{Dowling.2017} and \citet{Kraft.2022} do not consider continuous trading as relevant revenue stream of the intraday market.

On the continuous intraday market, market participants have to place their orders strategically since large volatilities driven by updates in the renewables forecast characterize the electricity price \citep{Kiesel.2017, Kremer.2020}.
For continuous trading, traders have to decide on the order type and the optimal sequence of order placements in a given time frame \citep{Guo.2017}.
This continuous trading can be modeled with the order-placement problem.
The order-placement problem is typically modeled with computational-expensive approaches, such as machine learning \citep{Boukas.2021}, dynamic programming \citep{Aid.2016, Glas.2020, Finnah.2021} and rolling-horizon optimization \citep{Corinaldesi.2020}, and allows to derive optimized trading strategies on the continuous intraday market.
These approaches focus on integrating uncertainties from renewables during continuous trading.
However, the approaches neglect other short-term markets.

The order-placement problem determines optimized trading strategies in the continuous intraday market. 
However, for optimal commitments to other markets, only potential revenues from utilizing the available flexible capacity in the continuous intraday market need to be determined one day in advance.
To value the trading of demand response flexibility in a real-time electricity market \citet{Muthirayan.2021} propose using option-price theory.
Thus, option-price theory can estimate the value of continuous trading under uncertain prices \citep{Weber.2015}.
As an advantage, the revenue of intraday trading can efficiently be estimated based on the price volatility without explicitly forecasting the future price curve.

\subsection{Contribution of this work}
Despite the growing importance of the continuous intraday market, the reviewed contributions on multi-market participation either neglect the continuous intraday market or the intraday market stage is represented without considering its key characteristic of continuous trading. 

Here, we present a coordinated bidding strategy that explicitly considers continuous trading in the intraday market while retaining a high level of detail for the day-ahead and balancing-power market.
The coordinated bidding strategy is derived from a multi-stage stochastic optimization model.
The multi-stage stochastic optimization models the sequential decision-making process for the balancing-power market, day-ahead market, and continuous intraday market.
For the continuous intraday market, we propose to use option-price theory to approximate the revenues from trading one day ahead.
Thus, the option-price theory allows to integrate trading on the continuous intraday market in a multi-stage stochastic optimization providing decision support in the sequential bidding process.

By combining option-price theory and multi-stage stochastic optimization, our method covers interdependencies between the balancing-power market, day-ahead market, and trading on the continuous intraday market to answer the question: Where to market flexibility?

\section{Method for optimized bidding strategy in sequential electricity markets}
\label{sec:method}
The proposed method allocates the flexible capacity of a market participant to the balancing-power market, the day-ahead market, and the continuous intraday market.
Thus, the method determines an optimized bidding strategy while also optimizing the operation of the considered energy system.
We base the method on our conference publication \citep{Nolzen.2022}.
Here, the method is extended by balancing-power market participation, including both stages for the acceptance and the request of balancing power.

Section~\ref{subsec:trading_deadlines} introduces the sequential decision-making process of multi-market participation.
We exemplarily present the different trading deadlines and how tenders are placed for the German balancing-power market, day-ahead market, and continuous intraday market.
However, the method is applicable to other countries with similar sequential electricity markets. 

Coordinated bidding requires to consider the value of the opportunity to trade on the continuous intraday market.
Following the ideas to model demand response \citep{Muthirayan.2021} and the continuous intraday market \citep{Weber.2015, Kern.2019b} as an option, we present an approach to value trading flexible capacity in the continuous intraday market in section~\ref{subsec:option_value}.
Therein, the value is derived using option-price theory. 

Coordinated bidding in sequential decision making leads to the stochastic process of multi-market participation (section~\ref{subsec:stochastic_process}).
In a final step, we set up a multi-stage stochastic optimization (section~\ref{subsec:stochastic_optimization}).

\subsection{Trading deadlines in short-term electricity and balancing-power markets}
\label{subsec:trading_deadlines}

This section presents the considered market design consisting of a balancing-power market, day-ahead market, and continuous intraday market. 
Thereby, we study the setting of sequential electricity markets that is common in Europe.
However, we adjust our method to the German electricity market design as Europe`s largest electricity market.
Thus, we consider the German electricity market design featuring the most common aspects of European electricity markets without limiting the general applicability of the method. 

\begin{figure}
	\includegraphics[width=16cm]{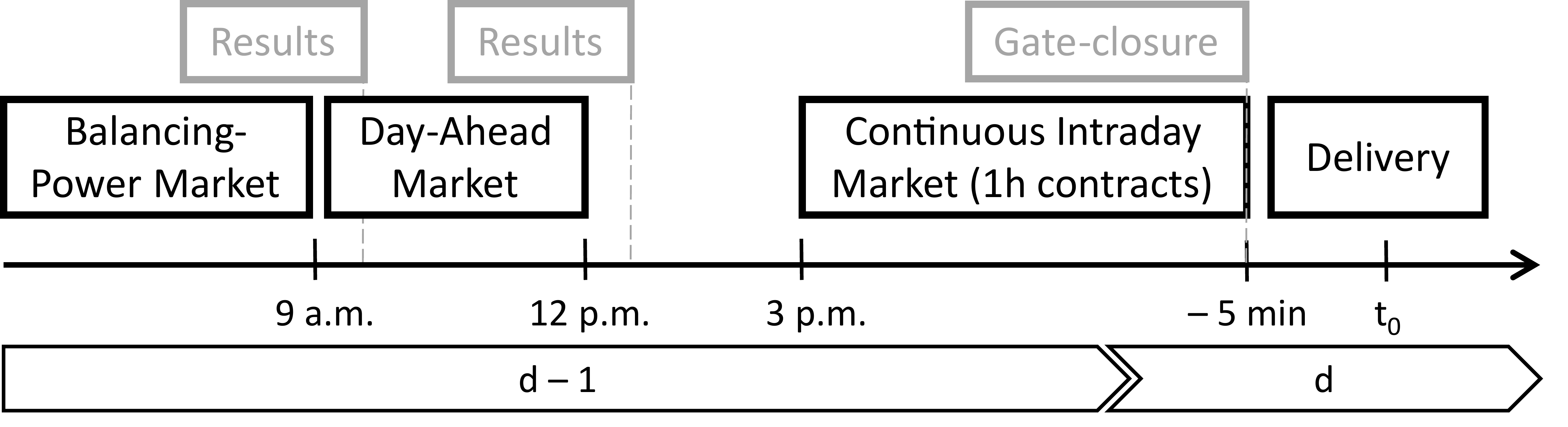}
 	\caption{Trading deadlines for sequential decision-making in balancing-power market, day-ahead market, and continuous intraday market oriented towards the German market design.}
	\label{fig:trading_deadlines}
\end{figure}

In the setting of sequential electricity markets, the balancing-power market, the day-ahead market, and the continuous intraday market are assigned different trading deadlines (Fig.~\ref{fig:trading_deadlines}).

Balancing power is traded on the balancing-power market, one day in advance (d--1).
We consider the market for automatic frequency restoration reserve (aFRR) with an activation time of 5 minutes.
On this balancing-power market, providers of positive and negative balancing power make a tender for four-hour time slices for the next day.
The tenders have to be submitted by 9 a.m. 
The auction results are published at 10 a.m.

Market participants submit a tender consisting of the amount of positive and/or negative balancing power and the respective capacity price bid and energy price bid.
The procurement of balancing power is organized as a daily two-stage pay-as-bid auction.
The first stage considers the capacity price bid that remunerates the provision of flexible capacity for balancing power.
The capacity price bids are sorted in ascending order.
Subsequently, the bids with the lowest capacity prices are accepted while the cut-off is made at the total demand of balancing power.
Among all market participants accepted in the first stage, the second stage considers the energy price bids.
The energy price bids remunerate the request of flexible capacity.
Again, the energy price bids are sorted in ascending order.
If balancing power is needed the next day, tenders with the lowest energy prices are requested first.
In both stages, market participants face a trade-off between prices and probabilities to maximize their revenues:
Low capacity prices lead to high acceptance probabilities but low capacity price revenues; 
Low energy prices lead to high request probabilities but low revenues per request.

\begin{figure}
	\includegraphics[width=16cm]{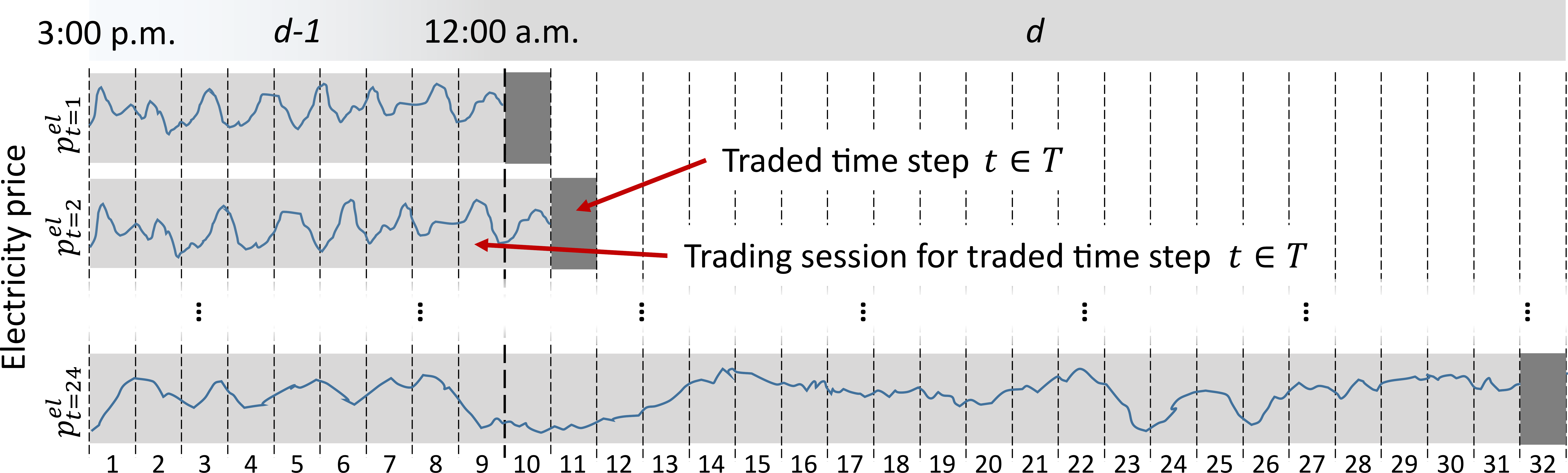} 
 	\caption{Trading sessions on the continuous intraday market for hourly contracts for the next day. The trading session for each of the contracts lasts from 9 hours to 32 hours.}
	\label{fig:trading_sessions_idm}
\end{figure}

On the day-ahead market, electricity is traded as hourly contracts for the next day.
The procurement is organized as a daily, double-sided blind auction. 
Buy-orders and sell-orders can be handed in until the order book closes at 12 p.m. 
Subsequently, aggregated electricity supply and demand curves are derived from the orders.
As a result, the market-clearing price is determined for each hour as the intersection of the demand and supply curves.
All market participants receive the market-clearing price which is the market equilibrium price.
The result of the day-ahead market auction is published at 1 p.m.

Trading on the continuous intraday market starts at 3 p.m.
Here, quarter-hourly contracts, half-hourly contracts, and hourly contracts are traded.
The contracts are continuously traded for the next day until shortly before delivery (gate-closure).
Thus, the trading sessions for each of the contracts last from 9 hours to 32 hours.
Fig.~\ref{fig:trading_sessions_idm} shows the trading session for hourly contracts.
Trading on the continuous intraday market follows the double-sided pay-as-bid principle.
The trade is executed at the bidding price, if a buy and sell order match.

In conclusion, we study a sequential market setting consisting of a two-stage pay-as-bid balancing-power market, a pay-as-clear day-ahead market, and a continuous intraday market, thus applying our method to a market design common in Europe for the method`s general applicability.  

\subsection{Approximating the revenues from the continuous intraday market}
\label{subsec:option_value}
For flexible capacity, the revenues in the continuous intraday market cannot be approximated with a single electricity price.
Additional revenues can be achieved with asset-backed trading as electricity is continuously sold and purchased depending on the price volatility \citep{Lohndorf.2021}.
Following \citet{Weber.2015}  and \citet{Muthirayan.2021}, we propose approximating the revenues from trading in the continuous intraday market with the option value.
This approximation has the advantage that the estimation is based on parameters available one day before delivery without explicitly forecasting the price curve.

We derive the option value with the multiperiod binomial model based on \citet{Cox.1979}.
The multiperiod binomial model models each trading session on the continuous intraday market in discrete time.
Therein, the option value is derived based on the stochastic price process during a trading session on the continuous intraday market and the marginal electricity costs for the flexible capacity.

Option price theory assumes a risk-neutral trader who aims to replicate the option value by the risk-neutral asset-backed trading strategy.
We refer to \ref{sec:abt_strategy} for an example of the risk-neutral asset-backed trading strategy since the explicit computation of this trading strategy is not necessary for the calculation of the option value, .

\vspace{3mm}
\textbf{Stochastic price process}
\vspace{2mm}

The option value requires the stochastic price process for each trading session as input parameter.
In a trading session, contracts are traded for the time step $t \in T$.
For each trading session, we discretize continuous trading with $N_t$ trading opportunities.
Therein, each trading session starts with the initial price level $S_t^\mathrm{ini}$.
For the initial price level, the day-ahead market price is considered as a suitable assumption \citep{Han.2022}.

Subsequently, the electricity price follows a stochastic price process.
Based on \citet{Alexander.2012}, we model the stochastic price process as an arithmetic Brownian motion that is able to capture negative electricity prices. 
Negative electricity prices have occurred more frequently in recent years, particularly during windy or sunny periods with a low residual load \citep{Halbrugge.2021}.
The arithmetic Brownian motion models the absolute price change by the price drift $\mu_t$ and the price volatility $\sigma_t$.
Thus, the arithmetic Brownian motion assumes that the absolute price deviation from the initial price level $S_t^\mathrm{ini}$ is normally distributed with $N(\mu_t,\sigma_t^2)$ at gate-closure.

The multiperiod binomial model discretizes the stochastic price process. 
Throughout a trading session, the price moves up $u_t$ or down $d_t$ at each trading opportunity.
For the underlying arithmetic Brownian motion, the up-movement $u_t$ and down-movement $d_t$ are determined as follows:
\begin{equation}
\label{eq: up_down_movement}
	\begin{aligned}
		u_t &= \mu_t \cdot \frac{1}{N_t} + \sigma_t\cdot\sqrt{\frac{1}{N_t}} \quad \forall t \in T,\\
		d_t &= \mu_t \cdot \frac{1}{N_t} - \sigma_t\cdot\sqrt{\frac{1}{N_t}} \quad \forall t \in T.
	\end{aligned}
\end{equation}

In Eq.~\eqref{eq: up_down_movement}, the price drift $\mu_t$ and the price volatility $\sigma_t$ are multiplied with the trading frequency $\frac{1}{N_t}$ to allocate the price drift $\mu_t$ and the price volatility $\sigma_t$ over the trading session.

The underlying arithmetic Brownian motion determines the last electricity price $S_{t,k}$ as the summation of up-movements and down-movements:
\begin{equation}
	\label{eq: arithmetic_brownian_motion}
	S_{t,k} = S_t^{\mathrm{ini}} + {k} \cdot {u_t} + {(N_t-k)} \cdot {d_t} \quad \forall t \in T, k\in \{ 0,1,…k,…,N_t \} .
\end{equation}
$N_{t}+1$ last prices $S_{t,k}$ are possible for each traded time step.

\vspace{3mm}
\textbf{Marginal electricity cost estimation}
\vspace{2mm}

The option value requires the marginal electricity costs $mc_t$ of the market participant as a second input parameter to estimate the revenues from risk-neutral asset-backed trading in the continuous intraday market (cf.~Fig.~\ref{fig:abt_strategy}). 
To determine the marginal electricity costs, we calculate the influence of changing electricity demand on the operational cost.
Therein, we do not allow for trading with the electricity and balancing-power markets because the flexibility needs to be provided by the energy system itself and not from the interaction with the markets.

We vary the electricity demand between zero and the maximum possible electricity output for the energy system.
For each electricity demand, we determine the cost-minimal operational schedule for the energy system by an operational optimization.
Afterward, a least-square linear regression is performed 
to obtain the operational cost as an affine function of the electricity demand.
The marginal electricity costs are approximated as the slope of the linear regression of the operational cost.
The linear approximation provides a good fit of the marginal electricity costs for the multi-energy system considered in section~\ref{sec:case_study}  ($R^2>0.97$).

The marginal electricity costs depend on the time step due to changing price parameters (e.g. gas price) and exogenous demands. 
For each time step $t \in T$, the procedure described above is repeated.

\newpage
\vspace{3mm}
\textbf{Option value as revenues from intraday market trading}
\vspace{2mm}

The option value is obtained based on the stochastic price process and the marginal electricity costs as input parameters.
The option value $opt_t^\mathrm{sell/buy}$ for positive (sell) and negative (buy) flexible capacity for each traded time step $t$ is derived from the multiperiod binomial model following \citet{Bjork.2009} with:
\begin{equation}
	\label{eq:option_price}
	opt_t^\mathrm{sell/buy}=\sum_{k=0}^{N_t} \underbrace{\binom{N_t}{k}}_{(a)} \cdot \underbrace{\left(\frac{-d_t}{u_t-d_t}\right)^k}_{(b)} \cdot \underbrace{\left(\frac{u_t}{u_t-d_t}\right)^{N_t-k}}_{(c)} \cdot \underbrace{\Phi^\mathrm{sell/buy}(S_{t,k})}_{(d)} \forall t\in T.
\end{equation}

Compared to \citet{Bjork.2009}, Eq.~\eqref{eq:option_price} is adapted to the continuous intraday market, setting the interest rate to zero due to the short-term nature of the market.
In Eq. \eqref{eq:option_price}, the binomial coefficient (a) denotes the absolute frequency that the last price $k$ is reached.
The factors (b) and (c) are the risk-neutral probabilities for an up-movement and down-movement.
The factor (d) evaluates the option value $\Phi^{sell/buy}(S_{t,k})$ at the end of a trading session for the sell option and the purchase option for each traded time step $t$ and last price $k$:
\begin{equation}
\Phi^\mathrm{sell}(S_{t,k})=
\begin{cases}
S_{t,k}-mc_t 	& \text{if} \quad S_{t,k} > mc_t \\
0 				& \text{if} \quad S_{t,k} \leq mc_t
\end{cases}	
\end{equation}
\begin{equation}
\Phi^\mathrm{buy}(S_{t,k})=
\begin{cases}
0 				& \text{if} \quad S_{t,k} \geq mc_t \\
mc_t - S_{t,k}	& \text{if} \quad S_{t,k} < mc_t.
\end{cases}
\end{equation}

In conclusion, the option value sums up and weights the revenues from all last prices $k$.
The option value $opt_t^\mathrm{sell/buy}$ represents the revenues from risk-neutral asset-backed trading for positive and negative flexible capacity.

As a result of risk-neutral asset-backed trading, either positive or negative flexible capacity is utilized.
Whether the respective flexible capacity is utilized, depends on the last price on the continuous intraday market compared to the marginal electricity costs. 
We use scenario probabilities $\pi_{t,\omega}^{\mathrm{ID}}$ that indicate the probability of utilizing positive and negative flexible capacity in the continuous intraday market.
For this purpose, we take the distribution of the last prices and evaluate its cumulative distribution function at the marginal electricity costs.
Then, the scenario probabilities $\pi_{t,\omega}^{\mathrm{ID}}$ are obtained from the cumulative distribution function indicating the probability that the last price is higher than the marginal electricity costs.

The option value $opt_t^\mathrm{sell/buy}$ and the scenario probabilities $\pi_{t,\omega}^{\mathrm{ID}}$ are the input parameters for the multi-stage stochastic optimization presented in the following sections~\ref{subsec:stochastic_process} and \ref{subsec:stochastic_optimization}.

\subsection{Stochastic process of multi-market participation}
\label{subsec:stochastic_process}

Participation in short-term electricity markets is a sequential decision-making process (cf. section~\ref{subsec:trading_deadlines}).
The sequential decision-making process is modeled as a stochastic process that considers all possible outcomes during market participation and operation.
Therein, we divide the stochastic process into four stages.
At each stage, either new information is revealed, or decisions are made. 

\begin{figure}
    \centering
	\includegraphics[width=8.2cm]{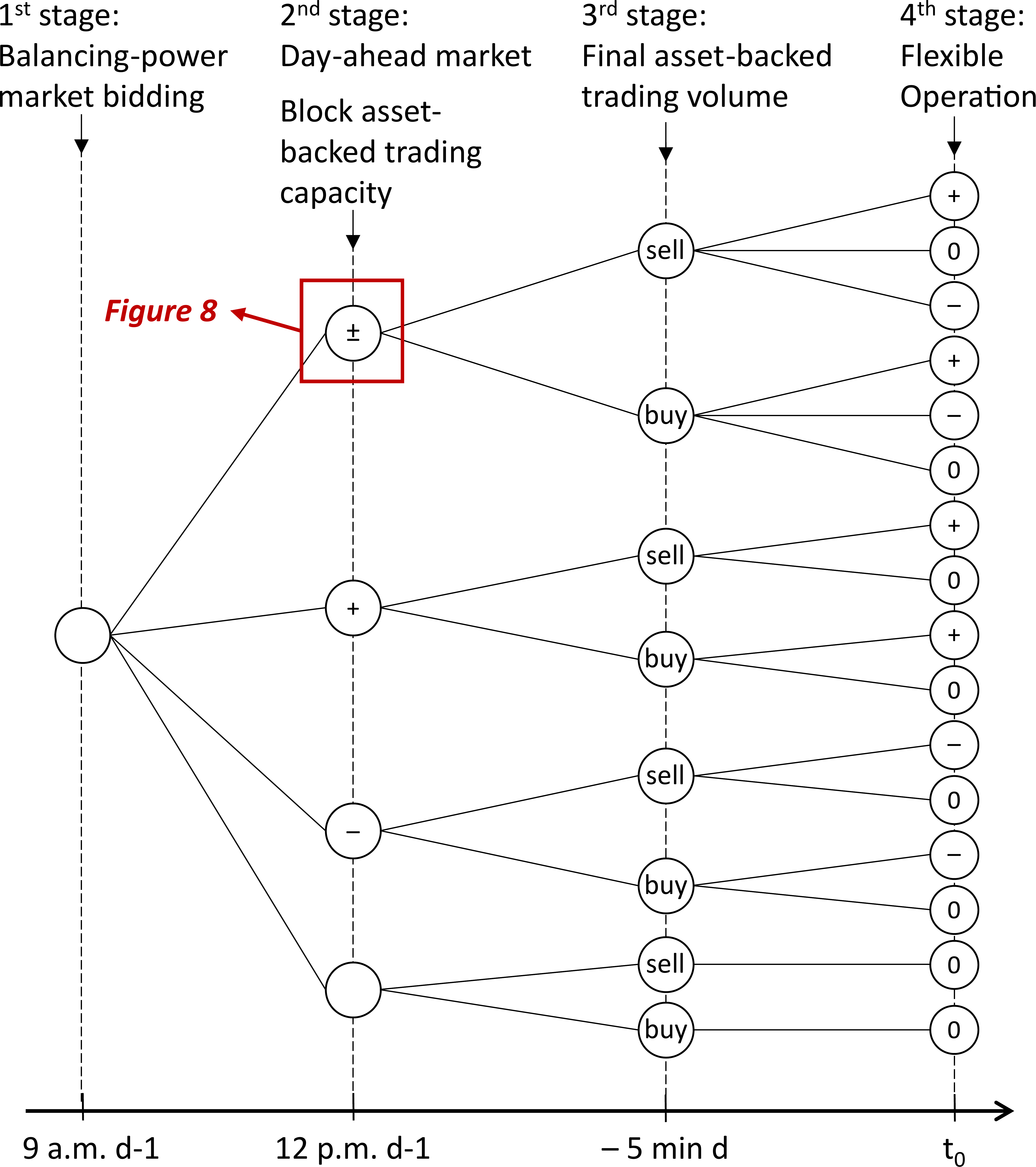}
 	\caption{Stochastic process of the multi-stage stochastic optimization model to obtain the optimal bidding strategy for a flexible market participant participating in the balancing-power market, day-ahead market, and continuous intraday market for each time step. Each stage reveals new information or makes a decision regarding market participation and operation.}
	\label{fig:stochastic_program}
\end{figure}

The resulting scenario tree comprises four stages and 16 scenarios $\omega \in \Omega$ for each time step $t \in T$ (Fig.~\ref{fig:stochastic_program}).
At the first stage (9 a.m. d--1), a tender is submitted for the balancing-power market.
The tender consists of a price combination $c$ with the capacity price ${cp}^{+/-}_{t,c}$ and energy price ${ep}^{+/-}_{t,c}$ and the amount of positive and negative balancing power ${BP}^{+/-}_t$.
After the results have been revealed at 10 a.m., the second stage (12 p.m. d--1) decides on how much electricity to sell ${DA}^{\mathrm{sell}}_{t,\omega}$, or buy ${DA}^{\mathrm{buy}}_{t,\omega}$ at the day-ahead market for the next day.
Electricity is purchased and sold at deterministic prices $p^{\mathrm{DA, buy}}_t$ and $p^{\mathrm{DA, sell}}_t$ since the energy system acts as a price taker.
The results on the delivery promises are revealed around 1 p.m.
Additionally, the second stage decides how much flexible positive capacity ${ID}^{\mathrm{sell}}_{t,\omega}$ (sell option) and negative capacity ${ID}^{\mathrm{buy}}_{t,\omega}$ (purchase option) are blocked for trading in the continuous intraday market.

The third stage (3 p.m. d--1 until gate-closure) considers asset-backed trading on the continuous intraday market with the amount of blocked flexible positive and negative capacity.
Depending on the marginal electricity costs and the electricity price development, two scenarios $\omega^{\mathrm{ID}} \in \{\text{sell}, \text{buy} \}$ arise for the final asset-backed trading volume.
Suppose the last price on the continuous intraday market is larger than the marginal electricity costs $(S_{t,k} > mc_t)$.
In this scenario, positive flexible capacity is entirely utilized, and electricity is sold on the continuous intraday market.
If the last price on the continuous intraday market is lower than marginal electricity costs $(S_{t,k} < mc_t)$, negative flexible capacity is entirely utilized, and electricity is purchased on the continuous intraday market.
At the fourth stage, the units $u \in U$ adapt the operation ${P}_{u,t,\omega,p}$  according to the request scenario that materializes.

\subsection{Modeling equations for the multi-stage stochastic optimization}
\label{subsec:stochastic_optimization}
This section presents the multi-stage stochastic optimization that models the stochastic process from the previous section \ref{subsec:stochastic_process}.
First, the objective function and the multi-market bidding problem are presented in section~\ref{subsec:objective_function}.
In section~\ref{subsec:energy_system_constraints}, the operational constraints of the underlying flexible multi-energy system are presented.
	
\subsubsection{Objective function and multi-market bidding problem}
\label{subsec:objective_function}
The objective function minimizes the expected operational cost ${OPEX}_t^{\mathrm{exp}}$ for all time steps~$t \in T$: 
\begin{equation}
	\min \sum_{t \in T}{OPEX}_t^{\mathrm{exp}}.
\end{equation}

In the balancing-power market, the acceptance probability and the request probability depend on the selected capacity price and energy price, respectively.
Hence, participation in the balancing-power market introduces endogenous uncertainty to the multi-stage stochastic optimization \citep{Hellemo.2018}.
As a result, the problem would be non-linear, usually requiring customized decomposition techniques, e.g., as in \citet{Schafer.2019}.

To avoid nonlinearities,
we extend the idea from \citet{Leenders.2020} instead: 
We introduce discrete capacity and energy prices with the set $C$.
The set $C$ contains each possible combination of capacity and energy prices for positive and negative balancing power.

We then model the scenario-dependent operational cost ${OPEX}_{t,c,\omega}^{\mathrm{scenario}}$ for all time steps~$t \in T$, price combinations~$c \in C$ and scenarios $\omega \in \Omega$ as
\begin{equation}
	\label{eq:OPEX-tcs}
	{OPEX}_{t,c,\omega}^{\mathrm{scenario}} = 
	C^{\mathrm{gas}}_{t,\omega}
	-R^{\mathrm{BP}}_{t,c,\omega}
	-R^{\mathrm{DA}}_{t,\omega}
	-R^{\mathrm{ID}}_{t,\omega}, \quad \forall \, t \in T, \, c \in C, \, \, \omega \in \Omega.
\end{equation}

The scenario-dependent operational costs comprise four terms: the cost of natural gas $C^{\mathrm{gas}}_{t,\omega}$, the revenues from balancing-power market participation $R^\mathrm{BP}_{t,c,\omega}$, the revenues from day-ahead market participation $R^\mathrm{DA}_{t,\omega}$, and the revenues from intraday market participation $R^{\mathrm{ID}}_{t, \omega}$.

The scenario-dependent operational costs ${OPEX}_{t,c,\omega}^{\mathrm{scenario}}$ are weighted with the scenario probability $\pi_{t,c,\omega}$.
Thereby, we obtain the operational costs ${OPEX}_{t,c}^{\mathrm{price~comb}}$ for each time step~$t \in T$, and price combination~$c \in C$ using:
\begin{equation}
	\label{eq:OPEX-tc}
	{OPEX}_{t,c}^{\mathrm{price~comb}} = \sum_{\omega \in \Omega} {OPEX}_{t,c,\omega}^{\mathrm{scenario}} \cdot \pi_{t,c,\omega}, \quad \forall \, t \in T, \, c \in C. 
\end{equation}
The scenario-dependent operational costs ${OPEX}_{t,c,\omega}^{\mathrm{scenario}}$ for each scenario $\omega \in \Omega$ are weighted with the scenario probability:
\begin{equation}
	\pi_{t,c,\omega} = \pi^{\mathrm{BP,cp}}_{t,c,\omega} \cdot \pi^{\mathrm{ID}}_{t,\omega} \cdot \pi^{\mathrm{BP,ep}}_{t,c,\omega}, \quad \forall \, t \in T, \, c \in C, \, \omega \in \Omega.
\end{equation}

The scenario probability $\pi_{t,c,\omega}$ is derived as product of the acceptance probability in the balancing-power market $\pi^{\mathrm{BP,cp}}_{t,c,\omega}$, the probability of positive or negative flexible capacity being utilized in the intraday market $\pi^{\mathrm{ID}}_{t,\omega}$ (cf. section \ref{subsec:option_value}), and the request probability of balancing power $\pi^{\mathrm{BP,ep}}_{t,c,\omega}$.
The request probability indicates the fraction of time that the energy system is requested within each time step $t$. 
Hence, the request probability is zero for all scenarios where the capacity price bid is not accepted.

In the bidding problem, one price combination $c$ is ultimately selected.
This bidding decision is modeled with the binary decision variable $\lambda^{\mathrm{BP}}_{t,c}$.
The variable $\lambda^{\mathrm{BP}}_{t,c}$ indicates if a price combination is selected ($\lambda^{\mathrm{BP}}_{t,c}=1$) or not ($\lambda^{\mathrm{BP}}_{t,c}=0$).
For each time step, one price combination is selected with:
\begin{equation}
	\sum_{c \in C} \lambda^{\mathrm{BP}}_{t,c} = 1, \quad \forall \, t \in T.
\end{equation}

The price combination with the lowest operational costs is identified with a Big-M formulation with $\mathrm{M}^{\mathrm{BP}}$ being a sufficiently large number:
\begin{equation}
\label{eq:big_m}
	{OPEX}_t^{\mathrm{exp}} \geq {OPEX}_{t,c}^{\mathrm{price~comb}}-\mathrm{M}^{\mathrm{BP}} \cdot (1-\lambda^{\mathrm{BP}}_{t,c}), \quad \forall \, t \in T, \, c \in C.
\end{equation}

As we consider a cost-minimization problem, inequality constraint Eq.~\eqref{eq:big_m} sets the expected operational costs $OPEX_t^{\mathrm{exp}}$ to the operational costs $OPEX_{t,c}^{\mathrm{price~comb}}$ with the price combination that has the minimal operational cost.

\vspace{3mm}
\textbf{The cost of natural gas}
\vspace{2mm}

The cost of natural gas $C^{\mathrm{gas}}_{t,\omega}$ for each time step~$t \in T$ and scenario~$\omega \in \Omega$ are derived by multiplying the required amount of natural gas ${BUY}_{t,\omega,\mathrm{gas}}$ and the natural gas price $p^{\mathrm{gas}}_t$
\begin{equation}
	C^{\mathrm{gas}}_{t,\omega} = {BUY}_{t,\omega,\mathrm{gas}} \cdot p^{\mathrm{gas}}_t, \quad \forall \, t \in T, \, \omega \in \Omega.
\end{equation}

\newpage
\vspace{3mm}
\textbf{The revenues from balancing-power market participation}
\vspace{2mm}

The revenues in the balancing-power market $R^\mathrm{BP}_{t,c,\omega}$ are modeled for each time step~$t \in T$, price combination~$c \in C$, and scenario~$\omega \in \Omega$ with:
\begin{equation}
\begin{split}
	R^\mathrm{BP}_{t,c,\omega}=&{BP}^+_t \cdot (s^{\mathrm{cp}+}_{\omega} \cdot {cp}^+_{t,c} + s^{\mathrm{ep}+}_{\omega} \cdot {ep}^+_{t,c}) \\
	+& {BP}^-_t \cdot(s^{\mathrm{cp}-}_{\omega} \cdot {cp}^-_{t,c} + s^{\mathrm{ep}-}_{\omega} \cdot {ep}^-_{t,c}), \quad \forall \, t \in T, \, \omega \in \Omega.
\end{split}
\end{equation}

The offered amount of positive balancing power ${BP}^+_t$ and negative balancing power ${BP}^-_t$ is compensated with the capacity price ${cp}^{+/-}_{t,c}$ and energy price ${ep}^{+/-}_{t,c}$ depending on the scenario $\omega \in \Omega$.
The binary parameters $s^{\mathrm{cp}+}_{\omega}$, $s^{\mathrm{cp}-}_{\omega}$, $s^{\mathrm{ep}+}_{\omega}$, $s^{\mathrm{ep}-}_{\omega}$ are equal to 1 if compensation is awarded in the respective scenario $\omega \in \Omega$.
The binary parameters $s^{\mathrm{cp}+}_{\omega}$ and $s^{\mathrm{cp}-}_{\omega}$ indicate the acceptance in the balancing-power market for positive and negative balancing power.
The parameters $s^{\mathrm{cp}+}_{\omega}$, $s^{\mathrm{cp}-}_{\omega}$ are equal to 1 if the capacity price bid for positive and negative balancing power is accepted in scenario $\omega$, and 0 otherwise.
The binary parameters $s^{\mathrm{ep}+}_{\omega}$ and $s^{\mathrm{ep}-}_{\omega}$ indicate the request of positive and negative balancing power. 
The parameters $s^{\mathrm{ep}+}_{\omega}$, $s^{\mathrm{ep}-}_{\omega}$ are equal to 1 if positive or negative balancing power is requested in scenario $\omega$, and 0 otherwise.
Hence, we assume a full request of balancing power (cf. \citep{Kumbartzky.2017, Bohlayer.2020b}).

Usually, balancing power is offered for several time steps.
For example, tenders are submitted for time slices of four hours in the German balancing-power market.
The price combination $c$ of capacity and energy price, as well as the offered amount of positive and negative balancing power have to be equal for each four-hour time slice.
Eq.~\eqref{eq:zeitscheibe_BP_pos} - \eqref{eq:zeitscheibe_lambda} ensure that the tender is the same for each four-hour time slice:
\begin{align}
{BP}^+_t&={BP}^+_{t+a}, \quad \forall \, t \in T: (t-1) \mid 4 = 0, \, a \in \{1,2,3\}, \label{eq:zeitscheibe_BP_pos}\\
{BP}^-_t&={BP}^-_{t+a}, \quad \forall \, t \in T: (t-1) \mid 4 = 0 , \, a \in \{1,2,3\}, \label{eq:zeitscheibe_BP_neg}\\
\lambda^{\mathrm{BP}}_{t,c}&=\lambda^{\mathrm{BP}}_{t+a,c}, \quad \forall \, t \in T: (t-1) \mid 4 = 0, \, a \in \{1,2,3\}. \label{eq:zeitscheibe_lambda}
\end{align}

\vspace{3mm}
\textbf{The revenues from day-ahead market participation}
\vspace{2mm}

The revenues from day-ahead market participation $R^\mathrm{DA}_{t,\omega}$ for each time step $t \in T$ and scenario $\omega \in \Omega$ equal the amount of sold electricity ${DA}^{\mathrm{sell}}_{t,\omega}$ and purchased electricity ${DA}^{\mathrm{buy}}_{t,\omega}$ multiplied with the day-ahead market price $p^{\mathrm{DA, sell}}_t$ for selling and $p^{\mathrm{DA, buy}}_t$ for purchasing electricity:
\begin{equation}
	R^\mathrm{DA}_{t,\omega}= 
	{DA}^{\mathrm{sell}}_{t,\omega} \cdot p^{\mathrm{DA, sell}}_t 
	+{DA}^{\mathrm{buy}}_{t,\omega} \cdot p^{\mathrm{DA, buy}}_t , \quad \forall \, t \in T, \, \omega \in \Omega.
\end{equation}

Electricity can either be sold or purchased at each time step $t \in T$ and scenario $\omega \in \Omega$.
We introduce the binary variable $\lambda^{\mathrm{DA}}_{t,\omega}$, which equals one if electricity is sold and equals 0 if electricity is purchased leading to
\begin{align}
	{DA}^{\mathrm{sell}}_{t,\omega} \leq \mathrm{M}^{\mathrm{DA}} \cdot \lambda^{\mathrm{DA}}_{t,\omega}, \quad \forall \, t \in T, \, \omega \in \Omega, \label{eq:only_sell}\\
	{DA}^{\mathrm{buy}}_{t,\omega} \leq \mathrm{M}^{\mathrm{DA}} \cdot (1 - \lambda^{\mathrm{DA}}_{t,\omega}), \quad \forall \, t \in T, \, \omega \in \Omega. \label{eq:only_buy}
\end{align}

Thereby, we restrict the purchase and selling of electricity at each time step $t \in T$ and scenario $\omega \in \Omega$ to the maximum possible amount of electricity $\mathrm{M}^{\mathrm{DA}}$ to be traded on the day-ahead market.

The sequential decision-making process in the markets is modeled with non-anticipativity constraints \citep{Birge.2011}.
Non-anticipativity constraints are used in stochastic programming to restrict the recourse of the multi-stage stochastic optimization.
Hence, decisions taken at previous stages cannot be changed at later stages.

The non-anticipativity constraints for the trading volumes on the day-ahead market are formulated as follows:
\begin{align}
	{DA}^{\mathrm{sell}}_{t,\omega} &= {DA}^{\mathrm{sell}}_{t,\omega'}, \quad \forall \, t \in T, \, 
	\omega,\omega' \in \Omega : (s^{\mathrm{cp}+}_{\omega},s^{\mathrm{cp}-}_{\omega})=(s^{\mathrm{cp}+}_{\omega'},s^{\mathrm{cp}-}_{\omega'}),
	 \label{eq:NAC_sell}\\
	{DA}^{\mathrm{buy}}_{t,\omega} &= {DA}^{\mathrm{buy}}_{t,\omega'}, \quad \forall \, t \in T, \,
	\omega,\omega' \in \Omega : (s^{\mathrm{cp}+}_{\omega},s^{\mathrm{cp}-}_{\omega})=(s^{\mathrm{cp}+}_{\omega'},s^{\mathrm{cp}-}_{\omega'}). 
	\label{eq:NAC_buy}
\end{align}

After the balancing-power market auction, Eq.~\eqref{eq:NAC_sell} and \eqref{eq:NAC_buy} ensure that the same trading volumes are submitted to the day-ahead market at each branch of the second stage.
If the binary parameter $s^{\mathrm{cp}+}_{\omega}$ and $s^{\mathrm{cp}-}_{\omega'}$ are equal in scenario $\omega$ and $\omega'$, Eq.~\eqref{eq:NAC_sell} and \eqref{eq:NAC_buy} equate the trading volumes on the day-ahead market.

\vspace{3mm}
\textbf{The revenues from intraday market participation}
\vspace{2mm}

The revenues from trading flexibility on the continuous intraday market $R^{\mathrm{ID}}_{t, \omega}$ are determined with the option value.
We assume that positive flexible capacity ${ID}^{\mathrm{sell}}_{t,\omega}$, and negative flexible capacity ${ID}^{\mathrm{buy}}_{t,\omega}$ can be monetized on the continuous intraday market.
Hence, the revenues in the intraday market $R^{\mathrm{ID}}_{t, \omega}$ are determined based on the amount of positive and negative flexible capacity ${ID}^{\mathrm{sell}}_{t,\omega}$ and ${ID}^{\mathrm{buy}}_{t,\omega}$ blocked for asset-backed trading, the respective option value ${opt}^{\mathrm{sell}}_t$ and ${opt}^{\mathrm{buy}}_t$, and the marginal electricity costs ${mc}_t$.

The option value is realized independently of the last price in the intraday market.
However, depending on the last price in the intraday market compared to the marginal electricity costs, two contrary scenarios $\omega^{\mathrm{ID}} \in \{\text{sell}, \text{buy} \}$ arise for the utilization of flexible capacity:
If the last price is higher than the marginal electricity costs, the multi-energy system supplies electricity to the grid to fulfill the delivery promise on the intraday market (sell).
This scenario utilizes the positive flexible capacity ${ID}^{\mathrm{sell}}_{t,\omega}$.
If the last electricity price is lower than the marginal electricity costs, the multi-energy system draws electricity from the grid to meet the purchase agreement on the intraday market (buy).
This scenario utilizes the negative flexible capacity ${ID}^{\mathrm{buy}}_{t,\omega}$.

Following, the revenues in the continuous intraday market are determined as:
\begin{equation}
	\label{eq:idm_revenues}
	\begin{split}
	R^{\mathrm{ID}}_{t, \omega} &= 
	{ID}^{\mathrm{sell}}_{t,\omega} \cdot ({opt}^{\mathrm{sell}}_t + s^{\mathrm{sell}}_{\omega} \cdot {mc}_t) \\
	&+ {ID}^{\mathrm{buy}}_{t,\omega} \cdot ({opt}^{\mathrm{buy}}_t - s^{\mathrm{buy}}_{\omega} \cdot {mc}_t) , \quad \forall \, t \in T, \, \omega \in \Omega.
	\end{split}
\end{equation}

We introduce the binary parameters $s^{\mathrm{sell}}_{\omega}$ and $s^{\mathrm{buy}}_{\omega}$ in Eq.~\eqref{eq:idm_revenues}.
The binary parameters model the scenario-dependent remuneration on the continuous intraday market and the physical electricity exchange with the continuous intraday market (cf. Eq.~\eqref{eq:electricity-balance}).

If the sell scenario materializes, the positive flexible capacity ${ID}^{\mathrm{sell}}_{t,\omega}$ is utilized.
The multi-energy system delivers electricity to the intraday market and receives the option value ${opt}^{\mathrm{sell}}_t$.
Since additional operational costs arise endogenously, the binary parameter $s^{\mathrm{sell}}_{\omega}$ is set to one.
We compensate the additional operational costs by the marginal electricity costs ${mc}_t$, as the operational costs from participating in the continuous intraday market are considered in the option value and in the stochastic optimization problem (cf. section~\ref{subsec:energy_system_constraints})
In the stochastic optimization problem, we thus avoid the operational costs being taken into account twice.

However, no electricity is purchased on the intraday market $(s^{\mathrm{buy}}_{\omega}=0)$ in the sell scenario, while the negative flexible capacity ${ID}^{\mathrm{buy}}_{t,\omega}$ is remunerated with the option value ${opt}_t^{\mathrm{buy}}$.  

If the buy scenario materializes, the negative flexible capacity ${ID}^{\mathrm{buy}}_{t,\omega}$ is utilized.
Electricity is purchased from the intraday market.
As the operational costs reduce endogenously, the binary parameter $s^{\mathrm{buy}}_{\omega}$ is set to one.
We subtract the saved operational costs by the marginal electricity costs ${mc}_t$.
In this scenario, the electricity delivery to the intraday market is not beneficial $(s^{\mathrm{sell}}_{\omega}=0)$.
However, the positive flexible capacity ${ID}^{\mathrm{sell}}_{t,\omega}$ is still remunerated with the option value ${opt}_t^{\mathrm{sell}}$.

Similar to Eq.~\eqref{eq:NAC_sell} and \eqref{eq:NAC_buy} for the day-ahead market, we introduce non-anticipativity constraints for the positive flexible capacity and negative flexible capacity blocked for trading on the continuous intraday market:
\begin{align}
	{ID}^{\mathrm{sell}}_{t,\omega} &= {ID}^{\mathrm{sell}}_{t,\omega'}, \quad \forall \, t \in T, \, 
	\omega,\omega' \in \Omega : (s^{\mathrm{cp}+}_{\omega},s^{\mathrm{cp}-}_{\omega})=(s^{\mathrm{cp}+}_{\omega'},s^{\mathrm{cp}-}_{\omega'}),
	 \label{eq:NAC_id_sell}\\
	{ID}^{\mathrm{buy}}_{t,\omega} &= {ID}^{\mathrm{buy}}_{t,\omega'} \quad \forall \, t \in T, \,
	\omega,\omega' \in \Omega : (s^{\mathrm{cp}+}_{\omega},s^{\mathrm{cp}-}_{\omega})=(s^{\mathrm{cp}+}_{\omega'},s^{\mathrm{cp}-}_{\omega'}). 
	\label{eq:NAC_id_buy}
\end{align}

\subsubsection{Energy system constraints}
\label{subsec:energy_system_constraints}

We exemplary show the model of the multi-energy system that we consider in the case study in section~\ref{sec:case_study}. 
However, the method applies to other flexible electricity-based systems, such as energy-intense processes or virtual power plants.

The multi-energy system consists of several units $u \in U$, which supply the demands of products $p \in P$ in time step $t \in T$ and scenario $\omega \in \Omega$.
We assume a flexible multi-energy system by neglecting constraints to restrict the operation (e.g., ramping constraints, minimum up-, and down-times, etc.).
Yet, we consider minimum part-load and load-dependent efficiencies of the units.
Since the multi-energy system is based on \citet{Baumgartner.2020}, we refer to this publication for the detailed model.

In this subsection, we solely present the electricity balance Eq.~\eqref{eq:electricity-balance} for time step $t$ and scenario~$\omega$ with:
\begin{multline}
\label{eq:electricity-balance}
	d_{\mathrm{el},t}
	- \sum_{u \in U} {P}_{u,t,\omega,\mathrm{el}}
	+ {DA}^{\mathrm{sell}}_{t,\omega} 
	+ s^{\mathrm{sell}}_\omega \cdot {ID}^{\mathrm{sell}}_{t,\omega}
	+ s^{\mathrm{ep}+}_\omega \cdot {BP}^+_t = \\
 	+ {DA}^{\mathrm{buy}}_{t,\omega} 
 	+ s^{\mathrm{buy}}_\omega \cdot {ID}^{\mathrm{buy}}_{t,\omega}
 	+ s^{\mathrm{ep}-}_\omega \cdot {BP}^-_t,
 	\quad \forall \, t \in T, \, \omega \in \Omega.
\end{multline}

In each time step $t$ and scenario $\omega$, the electricity demand comprises of the model exogenous demand $d_{\mathrm{el},t}$ and model endogenous demand and supply.
The variable ${P}_{u,t,\omega,\mathrm{el}}$ models the endogenous demand of electrically driven units and the supply of electricity generating units.

Additionally, the electricity balance considers market participation:
The electricity supply of the multi-energy system increases if electricity is sold on the day-ahead market ${DA}^{\mathrm{sell}}_{t,\omega}$, the intraday market ${ID}^{\mathrm{sell}}_{t,\omega}$, and if positive balancing power ${BP}^+_t $ is requested.
The electricity supply decreases if electricity is procured on the day-ahead market ${DA}^{\mathrm{buy}}_{t,\omega}$, the intraday market ${ID}^{\mathrm{buy}}_{t,\omega}$, and if negative balancing power ${BP}^-_t$ is requested.
Therein, the binary parameters $s^{\mathrm{ep}+}_{\omega}$, $s^{\mathrm{ep}-}_{\omega}$, $s^{\mathrm{sell}}_{\omega}$, $s^{\mathrm{buy}}_{\omega}$ are used to formulate the electricity balance for the different scenarios $\omega \in \Omega$.

Overall, the energy system constraints ensure the operational feasibility of the multi-energy system for time step $t \in T$ and scenario $\omega \in \Omega$ while providing flexibility.

\section{Case study for electricity market participation of a multi-energy system}
\label{sec:case_study}
\begin{figure}
	\includegraphics[width=16cm]{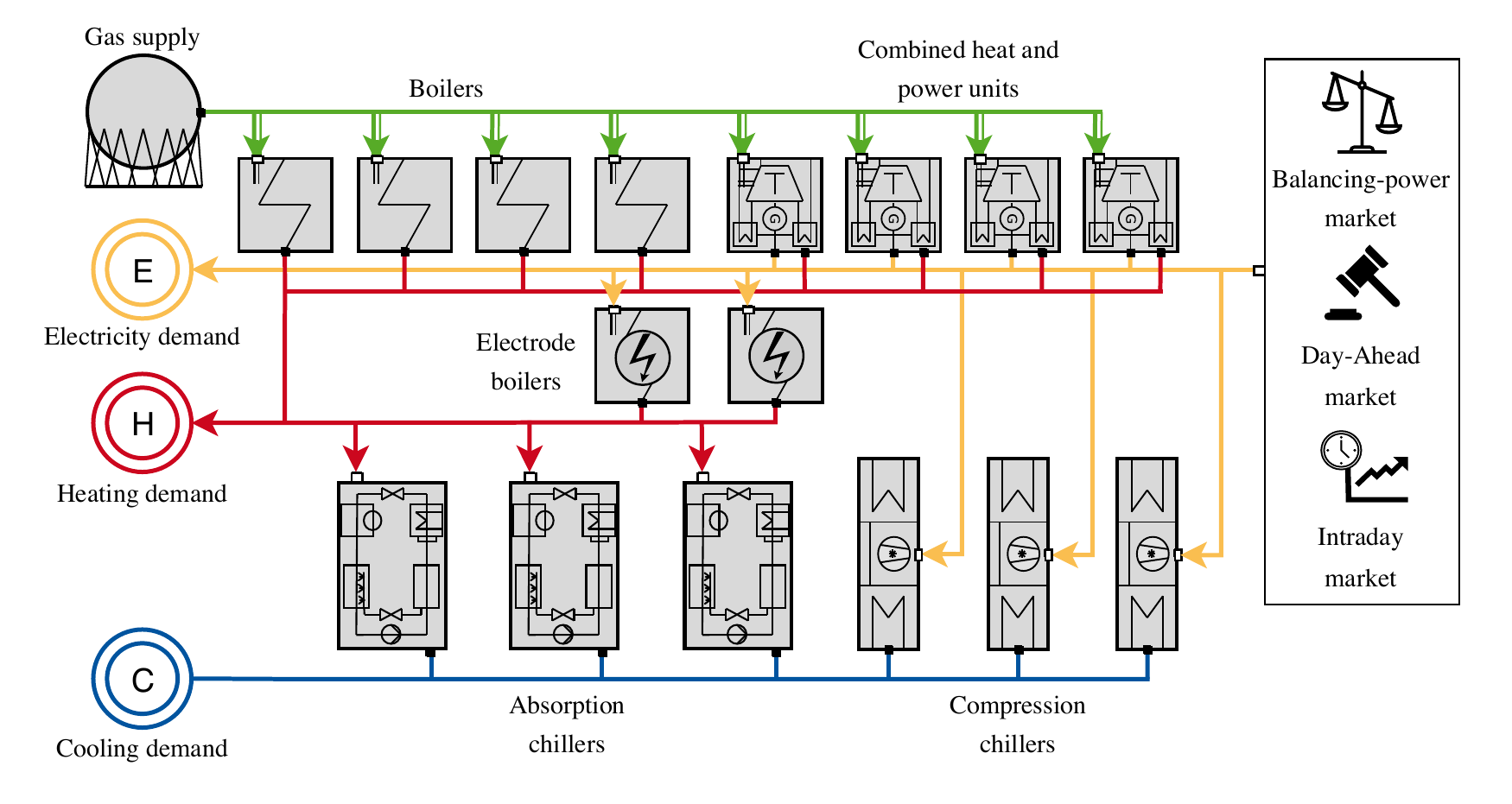}
 	\caption{Overview of the multi-energy system: absorption chillers, compression chillers, boilers, electrode boilers, and combined heat and power units cover hourly varying  electricity, heating, and cooling demands. The multi-energy system participates in the balancing-power market, day-ahead market, and continuous intraday market.}
	\label{fig:mes}
\end{figure}

The method proposed in section \ref{sec:method} is applied to a multi-energy system shown in Fig.~\ref{fig:mes} based on \citet{Baumgartner.2020}.
The multi-energy system provides flexibility by adjusting the operation of the units (Table~\ref{tab:components}) while supplying an industrial park with hourly varying electricity, heating, and cooling demands.
The multi-energy system participates in the German balancing-power market, day-ahead market, and continuous intraday market.
Therein, we base our analysis on the German electricity market design as of August 2019.

This section~\ref{sec:case_study} is structured as follows:
Section~\ref{subsec:inputs_modeling} presents the case study, while we prepare the market input parameters for our analysis in section~\ref{subsec:market_inputs_modeling}.
Subsequently, we present the results for the multi-market participation of the multi-energy system in section~\ref{subsec:market_participation} and perform a sensitivity analysis for the volatility in the continuous intraday market in section~\ref{subsec:volatility}.

We derive daily optimizations with 24 hourly time steps, each, within the SecMOD MILP framework \citep{Reinert.2022}.
All optimizations are solved using the solver Gurobi 9.5.0.
The optimizations are successfully solved within a time limit of \SI{3600}{\second} and a gap of 1~\% on an AMD EPYC 7F52 (16 cores, 3.5 GHz) with 192 GB RAM.

\subsection{Case study input parameters}
\label{subsec:inputs_modeling}
This section presents the input parameters for the case study.
For the balancing-power market, day-ahead market, and continuous intraday market, we base our study on historical data from one year (August 2019 -- July 2020) since no major changes in the balancing-power market design occurred during this period.

\begin{table}
	\caption{Overview of data sources to model the input parameters for the case study.}
	\centering	
	\begin{tabular}{ll} 
	\toprule
	\textbf{parameter} & \textbf{data source} \\
	\midrule
	electricity demand		& \citet{Baumgartner.2020} \\
	heating demand			& \citet{Baumgartner.2020} \\
	cooling demand			& \citet{Baumgartner.2020} \\
	\midrule 
	gas price 				& \citet{EPEXSPOT.2021} \\
	day-ahead market price	& \citet{BundesnetzagenturSMARD.de.2021} \\
	ID$_1$ price			& \citet{EPEXSPOT.2021} \\
	\midrule 
	Capacity price balancing-power market& \citet{Regelleistung.net.2021}\\  
	Energy price balancing-power market			& \citet{Regelleistung.net.2021} \\
	\midrule	
	Wind forecast 			& \citet{ENTSOE.2021} \\
	Photovoltaic forecast	& \citet{ENTSOE.2021} \\
	Total load forecast		& \citet{ENTSOE.2021} \\
	\bottomrule 
	\end{tabular}
\label{tab:price-data-sources}
\end{table}	

Tab.~\ref{tab:price-data-sources} lists the different data sources for the input parameters.
The gas price and the day-ahead market price are assumed to be deterministic since the prices are comparatively well predictable.
We neglect taxes and levies for the day-ahead market price as well as the gas price.
Thus, we consider the same price for selling and purchasing electricity on the day-ahead market.

For the continuous intraday market, we consider the ID$_1$ price to determine the price volatility.
The ID$_1$ price is the volume-weighted average price of the last hour of a trading session.

For the balancing-power market, we assume a maximum amount of \SI{10}{\mega\watt} positive and negative balancing power that can be offered as integer bids.
We take the marginal capacity price from \citet{Regelleistung.net.2021}.
The marginal capacity price is the highest price accepted in the capacity price auction.
We construct the marginal energy price based on the approach presented in \citet{Loesch.27.06.201829.06.2018}.
The marginal energy price is the highest price the request of balancing power is rewarded with.

\subsection{Preparation of market input parameters}
\label{subsec:market_inputs_modeling}
We derive the market input parameters based on the historical data and assumptions presented in the previous section~\ref{subsec:inputs_modeling}.
Our method requires the following market input parameters: the price drift and the price volatility of the continuous intraday market, discrete capacity and energy prices with corresponding acceptance and request probabilities for the balancing-power market. 
These parameters are partly dependent on each other since the balancing-power market and the continuous intraday market are interdependent \citep{Spodniak.2020}.
Hence, for a meaningful analysis of multi-market participation, it is important to derive consistent forecasts with consistent prices among the three considered markets \citep{Kraft.2022}.

Our focus is on the method to integrate trading on the continuous intraday market into the multi-market bidding problem.
Forecasting electricity prices is a research field of its own \citep{Hong.2020}.
While various methods exist to forecast day-ahead market prices \citep{Lago.2021}, a recent and active research area is forecasting intraday market prices \citep{Narajewski.2020} and balancing-power market prices \citep{ChristopherJahns.2019, Kraft.2020}.
Thereby, we acknowledge that the development of more accurate forecasting methods is out of the scope of this publication.

Various studies reveal the effect of wind generation \citep{Ketterer.2014}, solar generation \citep{Andrade.2017, Hu.2021}, and electricity demand \citep{Kiesel.2017} on market prices.
To prepare the price input parameters for our analysis, we select the wind generation forecast, solar generation forecast, and residual load forecast as influencing factors.
For those factors, the literature reveals a systematic influence on intraday prices, and a publicly available day-ahead estimate exists from \citet{ENTSOE.2021}.

\begin{figure}
    \centering
	\includegraphics[width=14cm]{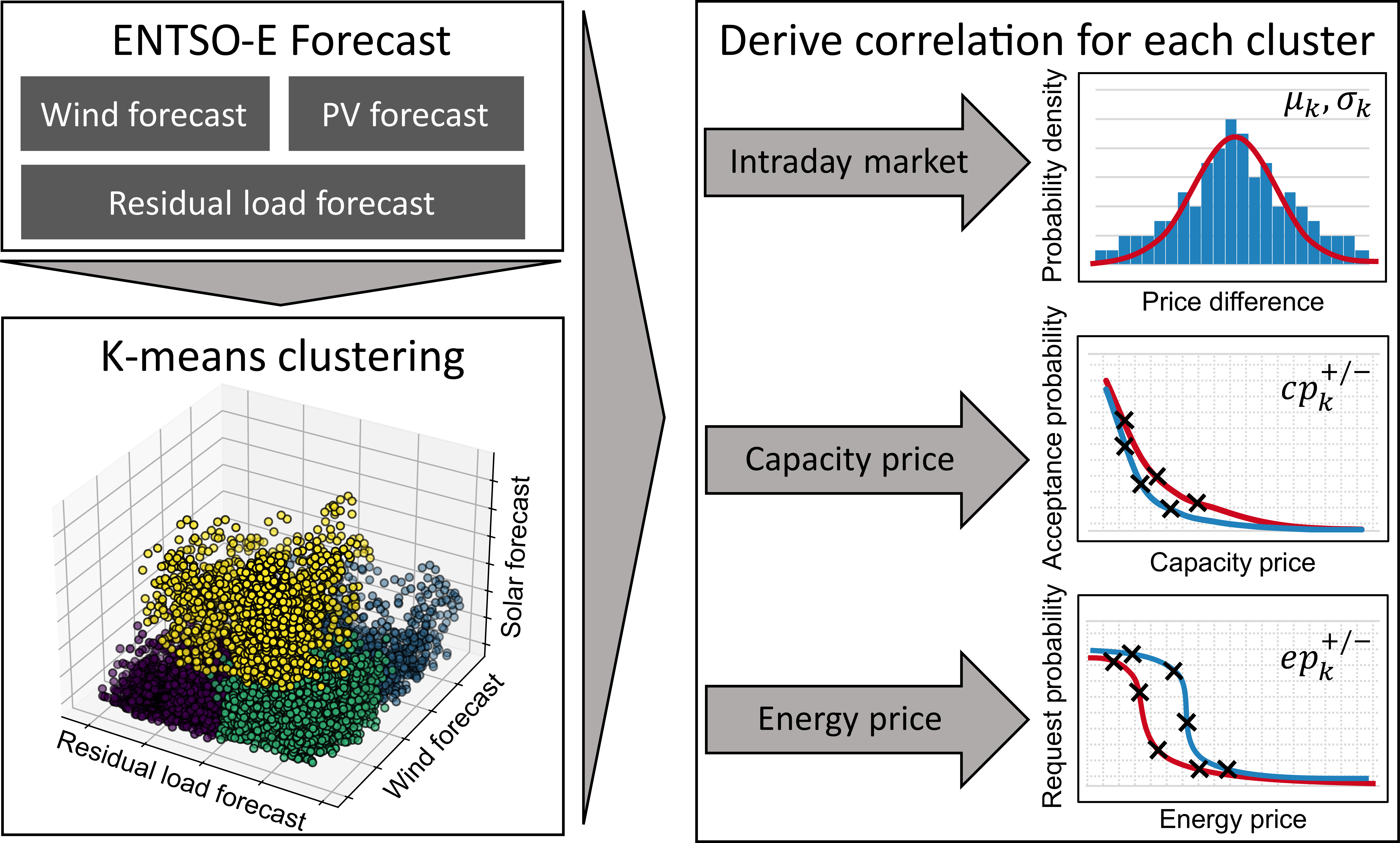}
 	\caption{Approach to derive the price input parameters for the balancing-power market, day-ahead market, and continuous intraday market. For the balancing-power market, we derive the cluster-dependent correlation and select discrete capacity prices $cp^{+/-}_k$ and energy prices $cp^{+/-}_k$. For the continuous intraday market, we obtain the cluster-dependent price drift $\mu_k$ and price volatility $\sigma_k$ to calculate the option value.}
	\label{fig:price_modeling_approach}
\end{figure}

The goal is now to obtain the market input parameters depending on the influencing factors to derive consistent price parameters among the considered markets (Fig.~\ref{fig:price_modeling_approach}).
Therefore, we first cluster historical data of the wind, photovoltaic and residual load forecast (Fig.~\ref{fig:price_modeling_approach} left) and then assign the market data to the respective clusters to obtain the market input parameters (Fig.~\ref{fig:price_modeling_approach} right).
Finally, the cluster-dependent market input parameters are assigned back to the respective time steps of the yearly time series.

On the day-ahead market and continuous intraday market, we consider hourly contracts.
Thus, we average the quarter-hourly forecast to hours as a first step.
Subsequently, we perform a k-means clustering of the historical data of the wind, photovoltaic and residual load forecast (Fig.~\ref{fig:price_modeling_approach} left).
For our analysis, we choose $|K|=4$, thus, dividing the data into 4 clusters.
Using more clusters, more differentiated conclusions could be drawn from the state of the electricity system regarding market prices.
However, more clusters decrease the data per cluster.
Our pre-analysis reveals $|K|=4$ as a good compromise between the number of clusters and data per cluster.
Now, the data in each cluster represents a different state of the electricity system resulting in different market input parameters.

On the balancing-power market, flexible capacity is provided for time slices of four hours.
Hence, we average the quarter-hourly forecast to 4 hours.
Subsequently, we also perform a k-means clustering with $|K|=4$. 	 

Following, we assign the market price data to each cluster (Fig.~\ref{fig:price_modeling_approach} right). For each cluster $k \in K$, we derive the price drift and price volatility in the continuous intraday market and price correlations for the balancing-power market as capacity price over the acceptance probability and the energy price over the request probability. 

For the continuous intraday market, we assign the price difference between the ID$_1$ price and the day-ahead market price to the clusters.
Subsequently, we fit a normal distribution to obtain the cluster-dependent price drift $\mu_k$ and price volatility $\sigma_k$.

For the balancing-power market, we derive the cluster-dependent acceptance probability $\pi^{\mathrm{BP}}_k ({cp}^{+/-})$ as a function of the capacity price.
For this purpose, we sort the capacity prices in descending order.
Thus, the acceptance probability is derived as:
\begin{equation}
	\pi^{\mathrm{BP}}_k ({cp}^{+/-}) = \frac{N_k({cp}^{+/-})}{N^{\mathrm{tot,cp}}_k}.
\end{equation} 
Therein, $N_k({cp}^{+/-})$ is the number of four-hour time slices with a marginal capacity price greater than the capacity price $cp^{+/-}$, while ${N^{\mathrm{tot,cp}}_k}$ is the total number of four-hour time slices in the considered cluster $k$.
The energy price is handled in the same way to derive the request probability $\pi^{\mathrm{BP}}_k ({ep}^{+/-})$ as a function of the energy price.

The multi-market optimization requires as inputs discrete capacity and energy prices with the respective acceptance and request probabilities.
The discrete prices are derived from the capacity price and energy price correlations.
From each capacity price correlation, we select three different capacity prices for both positive and negative balancing power.	
From each energy price correlation, we select four different energy prices for both positive and negative balancing power.
With this selection, we map the trade-off between capacity and energy prices, and their associated acceptance and request probabilities, respectively. 

More capacity and energy prices increase the number of binary variables in the bidding problem leading to higher computation times.
However, the number of selected capacity and energy prices still allows reasonable computational times of less than an hour.
The number of selected prices thus represents a trade-off between the computing time and the number of prices for the balancing-power market.

The demands of the considered multi-energy system are clustered with k-medoids.
Thereby, we obtain four typical weeks for our analysis in section~\ref{subsec:market_participation}:
calendar week 37 in 2019 (CW37-2019), calendar week 46 in 2019 (CW46-2019), calendar week 08 in 2020 (CW08-2020), and calendar week 18 in 2020 (CW18-2020).
For the four typical weeks, we compare the cases of participation in the day-ahead market (DA), day-ahead market and continuous intraday market (DA,ID), day-ahead market and balancing-power market (DA,BP), and day-ahead market, continuous intraday market and balancing-power market (DA,ID,BP).

\subsection{Results for the participation in sequential electricity markets}
\label{subsec:market_participation}
\begin{figure}
    \centering
	\includegraphics[width=9cm]{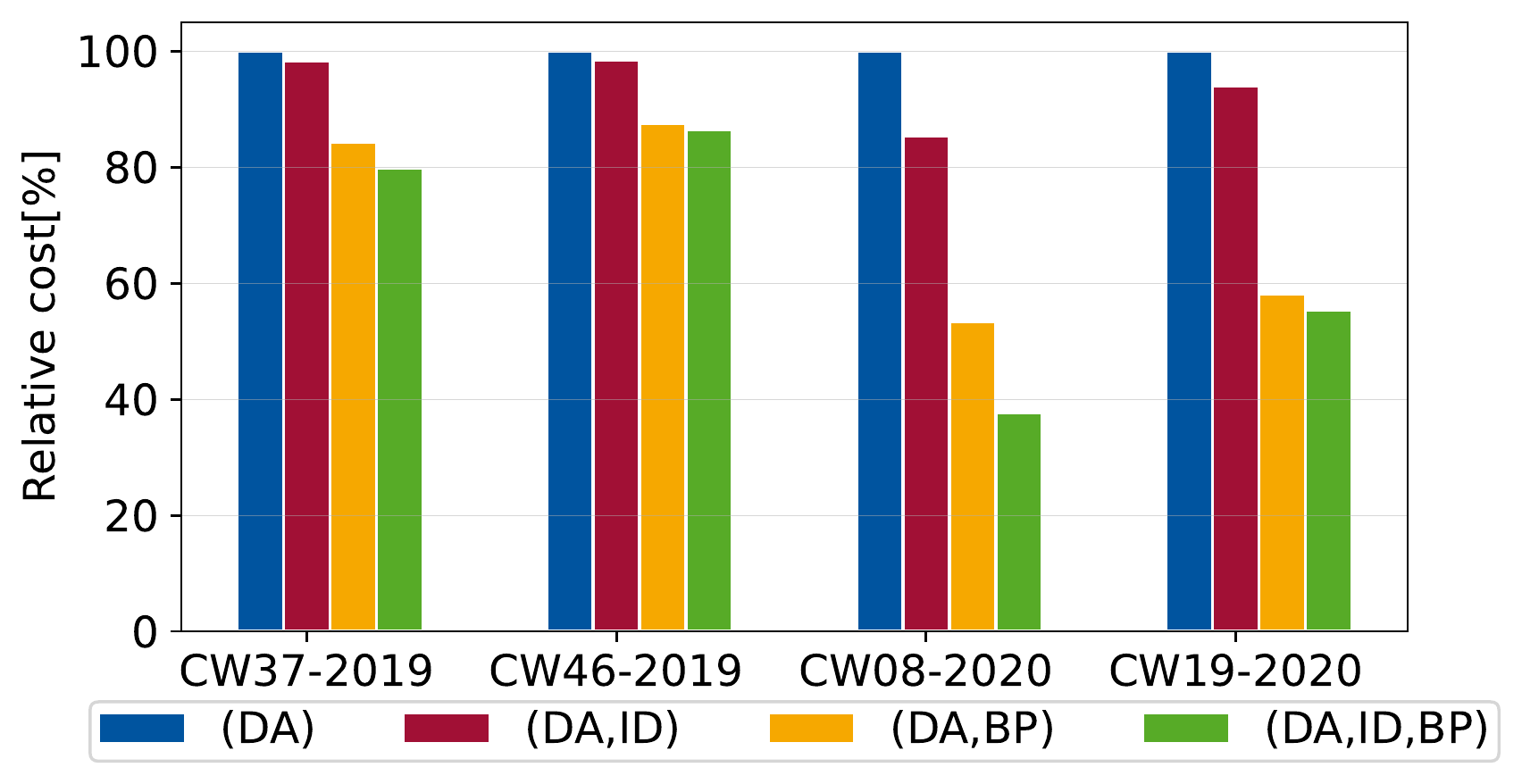}
 	\caption{Relative cost of the multi-energy system related to the costs in case (DA) for the four typical calendar weeks (CWs) in 2019/2020 when participating in the day-ahead market (DA), and in addition the continuous intraday market (DA,ID), or the balancing-power market (DA,BP), and in all three markets (DA,ID,BP).}
	\label{fig:relative_costs}
\end{figure}

The multi-energy system saves the most cost by simultaneous and coordinated participation in all three markets in case (DA,ID,BP) (Fig.~\ref{fig:relative_costs}).
In this case, the proposed method indicates savings of 14~\% to 62~\% in comparison to case (DA) for the four typical weeks.
With expanding participation in the day-ahead market to the continuous intraday market, costs decrease between 2~\% and 15~\% in case (DA,ID). 
Additional participation in the balancing-power market leads to savings between 12~\% and 47~\% for case (DA,BP).

In case (DA,ID,BP), the multi-energy system exploits the possibility to market flexibility on both the balancing-power market and the continuous intraday market.
The flexible capacity is preferably deployed in the balancing-power market, while the continuous intraday market is used as an opportunity in the event of an unsuccessful tender in the capacity price auction of the balancing-power market.
Therein, the day-ahead market is used to maximize the preferred direction of flexibility, i.e., electricity is sold on the day-ahead market to maximize negative flexibility and purchased to maximize positive flexibility. 
Thus, the method exploits the possibility of sequential bidding in the multi-stage stochastic optimization.

The continuous intraday market achieves the largest additional revenues in CW08-2020 (cf. case (DA) with (DA,ID) in Fig.~\ref{fig:relative_costs}).
In this week, the price volatility in the continuous intraday market is the highest due to many hours with high feed-in of wind generation.
We generally estimate the volatility rather conservatively based on the ID$_1$ price for hourly contracts (cf. section~\ref{subsec:market_inputs_modeling}).
Quarter-hourly contracts and utilization of shorter lead times result in higher volatilities with higher option values \citep{Han.2022}. 
Thus, the multi-energy system may achieve even higher revenues in the continuous intraday market.

\begin{figure}
    \centering
	\includegraphics[width=9cm]{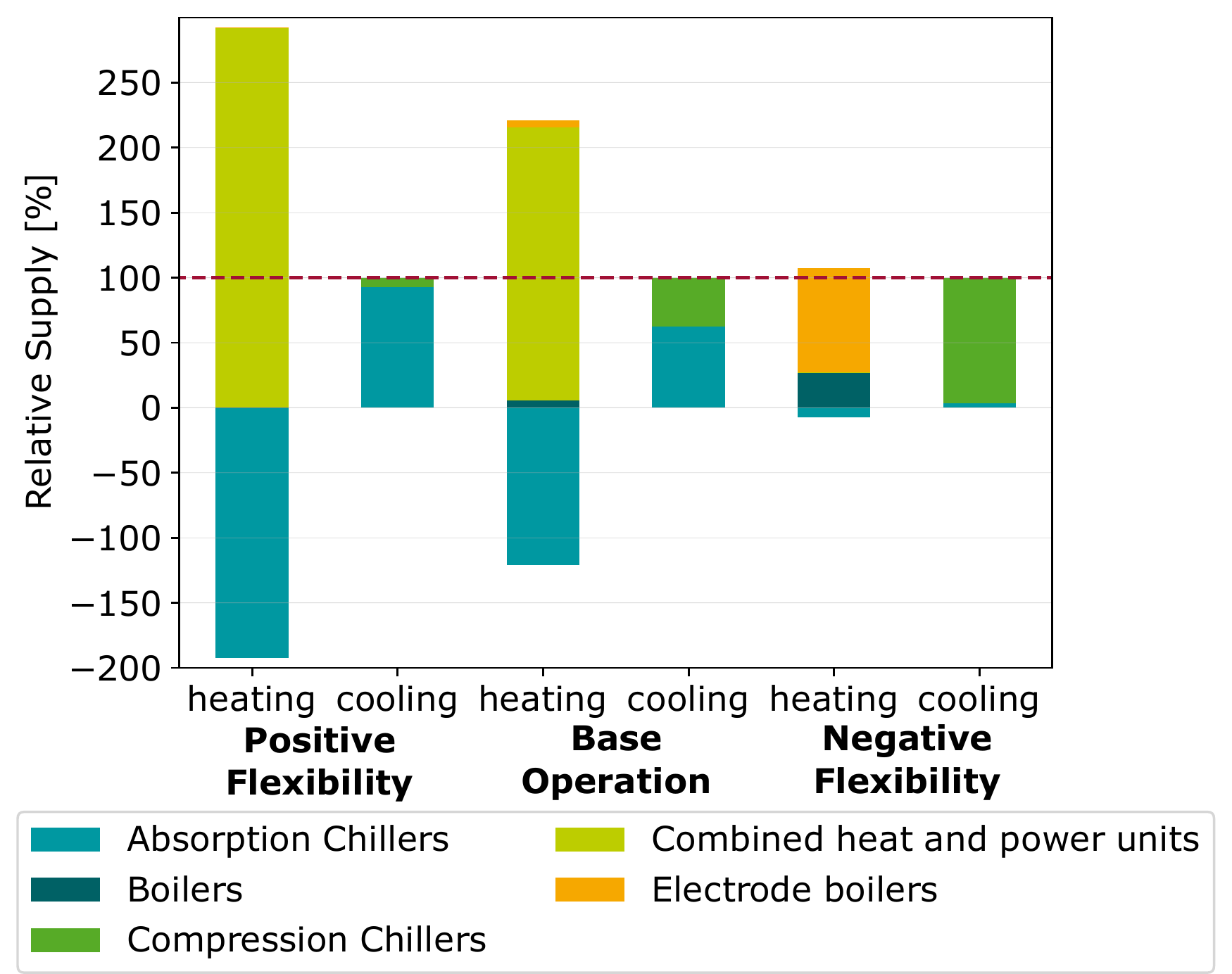}
 	\caption{Relative heating and cooling supply by the multi-energy system in calendar week 37 in 2019 (CW37-2019) for the utilization of the flexible capacity. Positive flexibility is provided for the request of positive balancing power and asset-backed trading of positive flexibility. Negative flexibility is provided for the request of negative balancing power and asset-backed trading of negative flexibility. Base operation refers to a scenario with no request of balancing power.}
	\label{fig:operation}
\end{figure}

The multi-energy system utilizes the flexible capacity by a complex and interconnected adjustment of its operation:
The electricity production is adjusted based on the market signals, followed by adjusting the operation of the multi-energy system for the heating and cooling production (Fig.~\ref{fig:operation}).
In the event that positive flexible capacity is utilized, the multi-energy system commonly increases electricity production by combined heat and power units.
Due to their co-generation, electricity production increases heat production, which is used to cover heating demands or cooling demands via absorption chillers.
In the event that negative flexible capacity is utilized, the multi-energy system either increases internal electricity use or decreases electricity production, as electricity is drawn from the grid.
Subsequently, electrode boilers mostly provide heat, and compression chillers provide cooling.

In conclusion, the results indicate significant savings when participating in all markets since the flexible capacity is utilized both in the balancing-power market and the continuous intraday market. 
The employed bidding strategies and operations are complex, justifying the use of the proposed method.

\subsection{Sensitivity analysis of intraday market volatility}
\label{subsec:volatility}

The previous analysis reveals that the utilization of flexible capacity is highly beneficial. 
However, the balancing-power market tends to be preferred for the typical weeks considered in this study.
To better understand market preferences, we analyze the sensitivities of the model with respect to the volatilities on the intraday market.
Higher volatilities on the intraday market add to the profitability of asset-backed trading and increase the option value of the intraday market in our model. 

\begin{figure}
    \centering
	\includegraphics[width=9cm]{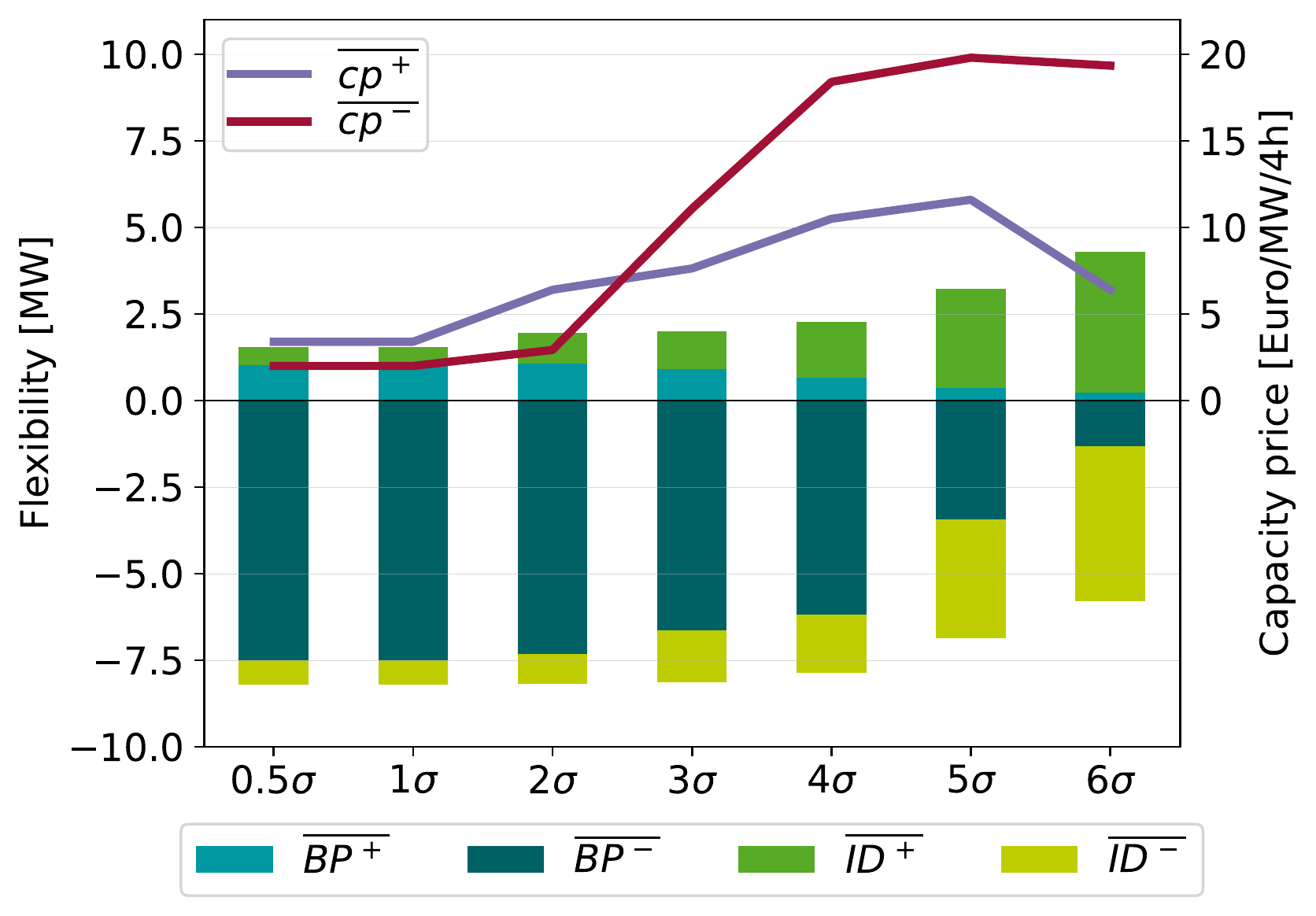}
 	\caption{Flexible capacity offered on the balancing-power market $\overline {BP^{+/-}}$ and continuous intraday market $\overline {ID^{+/-}}$ for CW37-2019 after the successful capacity price auction (cf. Fig.~\ref{fig:stochastic_program}) as a function of the intraday volatility $\sigma$. Additionally, the volume-weighted average capacity price $\overline {cp^{+/-}}$ is shown at which positive and negative balancing power is offered.}
	\label{fig:sensitivity_vola}
\end{figure}

We perform the sensitivity analysis for CW37-2019 and consider the case with participation in all markets (DA,ID,BP).
Fig.~\ref{fig:sensitivity_vola} compares the average amount of flexible capacity offered in the balancing-power market and continuous intraday market, assuming a successful capacity price auction for positive and negative balancing power.  
Thereby, we analyze both the balancing-power market bids of the first stage and the subsequent strategy on the continuous intraday market at the same time.

The sensitivity analysis highlights the importance of the intraday opportunity within a multi-market optimization (Fig.~\ref{fig:sensitivity_vola}):
Higher volatilities increase both the flexible capacity traded on the intraday market and the capacity prices in the balancing-power market, reducing the acceptance probabilities.
With increasing volatility the capacity prices are increased first ($2\sigma - 4\sigma$). 
Subsequently, the flexible capacity is shifted towards the continuous intraday market ($4\sigma - 6\sigma$).

We vary the volatility between a factor of 0.5 to 6 to reflect the opportunity for flexible capacity to utilize higher volatilities by trading quarter-hourly contracts and exploiting shorter lead times. 
With the base volatility ($1\sigma$), flexible capacity is mainly utilized in the balancing-power market (88~\%), while low capacity prices are offered. 
In this case, the average acceptance probability is 66~\% for positive and 48~\% for negative balancing power.
With increasing volatility, the offered capacity prices increase, and the flexible capacity is shifted towards the continuous intraday market. 
Consequently, the average acceptance probability reduces to 48~\% for positive and 10~\% for negative balancing power.
Thus, the provision of balancing power is less likely. 
In case $6\sigma$, almost the entire flexible capacity (85~\%) is instead allocated to the continuous intraday market.

In conclusion, the sensitivity analysis of the intraday volatilities shows that both the capacity prices and the flexible capacity are shifted to favor the continuous intraday market.
Thereby, the bidding strategy reflects the multi-energy system`s capabilities to exploit the volatility of the continuous intraday market.
To optimize the trade-off to utilize the flexible capacity, the presented method derives a coordinated bidding strategy among the markets.

\section{Conclusions}
\label{sec: conclusions}

A method is presented to integrate trading on the continuous intraday market in a multi-market optimization for the day-ahead, balancing-power, and continuous intraday market.
To approximate the revenues from continuous trading, we employ option-price theory. 
The option price is integrated into a multi-stage stochastic optimization that models multiple markets, including the energy system's operations.
As a result, the method is able to determine an optimal coordinated bidding strategy for all markets.
The method is not limited to energy systems as market participants.
Although we consider widely adopted electricity and balancing-power markets, the method can be adapted to different (future) market designs.

We apply our method to a multi-energy system participating in the German market as of August 2019.
The largest cost reductions can be achieved when participating in all three considered markets.
The method is able to utilize the flexible capacity of the multi-energy system.
A sensitivity analysis shows that larger volatilities in the intraday market support a shift towards increased participation in the continuous intraday market: Tenders on the balancing-power market increase capacity prices and thus reduce acceptance probabilities.
The analysis highlights the importance of considering the intraday opportunity.

In conclusion, our method derives a coordinated and optimized bidding strategy in short-term electricity and balancing-power markets by exploiting the flexibility of the market participant's energy system to achieve the largest cost reductions.

The performance of our method depends on sufficient and accurate price forecasts.
Since price forecasting is not in the scope of this paper, future research could focus on developing price forecasting methods that are consistent across the markets.  
In particular, the volatility on the continuous intraday market greatly influences market preferences.
Hence, forecasting methods for the continuous intraday market are desirable.
Additionally, we consider risk-neutral decision-making.
Since energy supply companies often act risk-averse, our approach could thus be extended by considering risk in the optimization model \citep{Kraft.2022,Germscheid.2022}.
Finally, our approach could consider storage to increase the flexibility of an energy system.

We hope that our method and the results motivate market participants to implement joint optimization in all markets. 
Joint market participation increases the flexibility of the power system and ultimately supports the energy transition.

\section*{Author contributions}
\textbf{Niklas Nolzen}: Conceptualization, Methodology, Software, Validation, Investigation, Data Curation, Writing - Original Draft, Visualization, Project administration. 
\textbf{Alissa Ganter}: Conceptualization, Methodology, Software, Investigation, Data Curation, Visualization.
\textbf{Nils Baumgärtner}: Conceptualization, Writing - Review \& Editing. 
\textbf{Ludger Leenders}: Conceptualization, Writing - Review \& Editing, Supervision, Visualization, Funding acquisition.
\textbf{André Bardow}: Conceptualization, Resources, Writing - Review \& Editing, Supervision, Funding acquisition.

\section*{Declaration of Competing Interest}
We have no conflict of interest.

\section*{Acknowledgements}
NN thanks the German Federal Ministry of Economic Affairs and Energy (ref. no.: 03EI1015A). AB and LL thank the Swiss Federal Office of Energy’s “SWEET” programme and the “PATHFNDR” consortium. The support is gratefully acknowledged. The authors thank Kang Qiu for her contribution to the market model and Laura Block for her contribution to the data preparation.

\newpage
\twocolumn
\section*{Nomenclature} \label{sec:nomenclature}

\noindent\textbf{Sets} \\
\noindent
\begin{tabularx}{\columnwidth}{lX}
$T$						& Time steps	\\
$\Omega$				& Scenarios		\\
$\Omega^{\mathrm{ID}}$ 	& Intraday market scenarios \\
$C$						& Capacity and energy price combinations \\
$U$						& Units						\\
$P$						& Products					\\
\end{tabularx} \\

\noindent\textbf{Variables} \\
\noindent
\begin{tabularx}{\columnwidth}{lX}
${OPEX}_t^{\mathrm{exp}}$					& Expected operational cost \\
${OPEX}_{t,c}^{\mathrm{price~comb}}$		& Operational cost per price combination \\
${OPEX}_{t,c,\omega}^{\mathrm{scenario}}$	& Operational cost per price combination and scenario \\
$C^{\mathrm{gas}}_{t,\omega}$	& Natural gas cost					\\
$R^{\mathrm{BP}}_{t,c,\omega}$	& Balancing-power market revenues	\\
$R^{\mathrm{DA}}_{t,\omega}$	& Day-ahead market revenues \\
$R^{\mathrm{ID}}_{t,\omega}$	& Intraday market revenues \\
$\lambda^{\mathrm{BP}}_{t,c}$	& Binary decision for price combination in the balancing-power market \\
$\lambda^{\mathrm{DA}}_{t,\omega}$	& Binary decision to buy/sell electricity in the day-ahead market \\
${BUY}_{t,\omega,p}$				& Amount of bought product (e.g. natural gas)			\\
${DA}^{\mathrm{sell}}_{t,\omega}$	& Amount of sold electricity in the day-ahead market	\\
${DA}^{\mathrm{buy}}_{t,\omega}$	& Amount of purchased electricity in day-ahead market 	\\
${ID}^{\mathrm{sell}}_{t,\omega}$	& Continuous intraday market trading volumes for positive flexible capacity \\
${ID}^{\mathrm{buy}}_{t,\omega}$	& Continuous intraday market trading volumes for negative flexible capacity \\
${BP}^+_t$				& Amount of positive balancing power\\
${BP}^-_t$				& Amount of negative balancing power\\
${P}_{u,t,\omega,p}$	& Product supply/demand of unit \\
\end{tabularx} \\

\newpage

\noindent\textbf{Parameters} \\
\noindent
\begin{tabularx}{\columnwidth}{lX}
$\pi_{t,c,\omega}$ 									& Scenario probability \\
$\pi^{\mathrm{BP,cp}}_{t,c,\omega}$ 				& Acceptance probability in the balancing-power market \\
$\pi^{\mathrm{ID}}_{t,\omega}$ 	& Probability of positive or negative flexible capacity being utilized in the intraday market  \\
$\pi^{\mathrm{BP,ep}}_{t,c,\omega}$					& Request probability in the balancing-power market \\
$\mathrm{M}^{\mathrm{BP}}$		& Big M (large number) \\
$\mathrm{M}^{\mathrm{DA}}$		& Big M (large number) \\
$p^{\mathrm{gas}}_t$			& Natural gas price \\
$p^{\mathrm{DA, sell}}_t$		& Day-ahead market selling price \\
$p^{\mathrm{DA, buy}}_t$		& Day-ahead market purchasing price \\
${cp}^+_{t,c}$					& Positive capacity price \\
${cp}^-_{t,c}$					& Negative capacity price \\
${ep}^+_{t,c}$					& Positive energy price \\
${ep}^-_{t,c}$					& Negative energy price \\
$s^{\mathrm{cp}+}_{\omega}$		& Binary parameter for acceptance of positive balancing power \\
$s^{\mathrm{cp}-}_{\omega}$		& Binary parameter for acceptance of negative balancing power \\
$s^{\mathrm{ep}+}_{\omega}$		& Binary parameter for request of positive balancing power \\
$s^{\mathrm{ep}-}_{\omega}$		& Binary parameter for request of negative balancing power \\
$s^{\mathrm{sell}}_{\omega}$	& Binary parameter to sell in the continuous intraday market \\
$s^{\mathrm{buy}}_{\omega}$		& Binary parameter to buy in the continuous intraday market \\
$d_{p,t}$						& Product demand \\
${opt}^{\mathrm{sell}}_t$		& Option price for positive flexible capacity \\
${opt}^{\mathrm{buy}}_t$		& Option price for negative flexible capacity \\
${mc}_t$						& Marginal electricity cost \\
$u_t$							& Up-movement \\
$d_t$							& Down-movement \\
$\mu_t$ 						& Price drift \\
$\sigma_t$						& Price volatility \\
$N_t$							& Number of trading opportunities \\
$S_t^\mathrm{ini}$				& Initial price-level				\\
$S_{t,k}$						& Last prices						\\
$\Phi^\mathrm{sell}$			& Option value for positive flexible capacity at the end of the trading session \\
$\Phi^\mathrm{buy}$				& Option value for negative flexible capacity at the end of the trading session \\
\end{tabularx} \\

\newpage 
\onecolumn
\appendix
\label{Appendix}
\section{Risk-neutral asset-backed trading strategy}
\label{sec:abt_strategy}

\begin{figure}[h!]
    \centering
	\includegraphics[width=9cm]{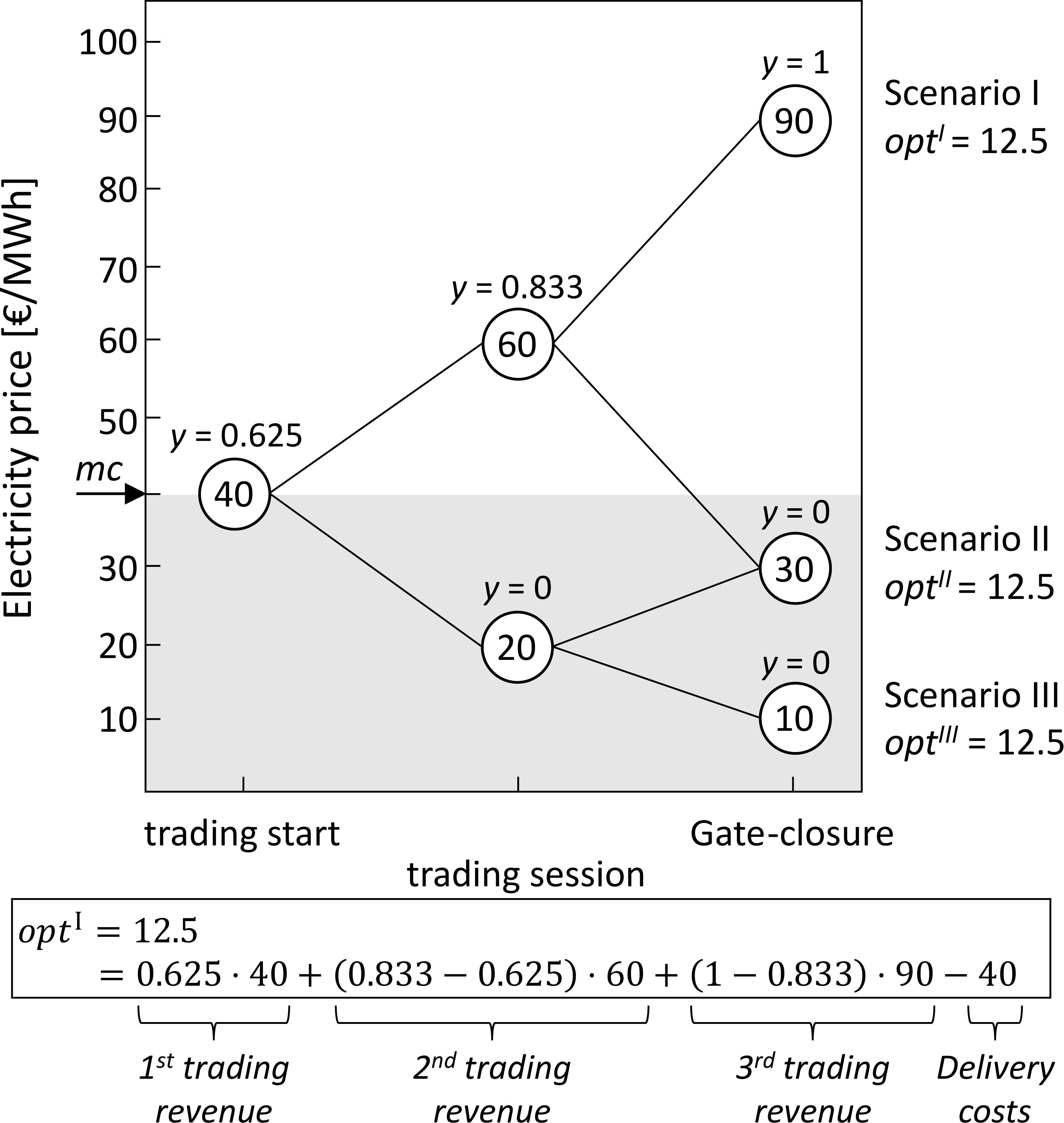}
 	\caption{Exemplary risk-neutral asset-backed trading strategy for positive flexible capacity in the multiperiod binomial model based on \citep{Cox.1979}. Each discrete price offers a trading opportunity at which the trader adjusts the share of sold flexible capacity $y$ to realize the option value $opt$ independent of scenarios I, II, and III.}
	\label{fig:abt_strategy}
\end{figure}

Fig.~\ref{fig:abt_strategy} exemplifies the risk-neutral asset-backed trading strategy in the multiperiod binomial model for positive flexible capacity.
For negative flexible capacity, the multiperiod binomial model works similarly but in the opposite direction.

The relative electricity output $y$ is the relative flexible capacity traded on the continuous intraday market.
At each trading opportunity, the trader adjusts the sold share of flexible capacity $y$.
The adjustment is based on the stochastic price process, and the marginal electricity costs $mc$.
Therein, the trader behaves risk-neutral to realize the estimated option value $opt$ within the trading session.
The option value is realized independent of the price scenario I, II, and III.
At the end of a trading session, the price volatility in the continuous intraday market is monetized to the option value.
The positive flexible capacity is either entirely sold ($y=1$) or not ($y=0$).

For further information regarding the explicit computation of the risk-neutral asset-backed trading strategy, we refer to \citet{Bjork.2009}.

\clearpage

\section{Installed units of the multi-energy system}

\begin{table}[!h]
	\caption[Units of the decentralized energy system]{Overview of the installed units of the multi-energy system, their electrical and thermal capacities, and their output. AC: absorption chillers, CC: compression chillers, B: boilers, EB: electrode boilers, and CHP: combined heat and power units.}
	\vspace{6pt}
	\centering	
	\begin{tabular}{llll} 
	\toprule
	\textbf{Component}	&	\textbf{Capacity} [\si{\mega\watt}\textsubscript{el}] &	\textbf{Capacity} [\si{\mega\watt}\textsubscript{th}]	&	\textbf{Output}\\
	\midrule 
	AC 1 		& &	3.5	&	Cooling\\
	AC 2		& &	2.5 	&		\\
	AC 3		& &	0.5 	&		\\
	\midrule 
	CC 1		& 0.8 &	4.5 	&	Cooling\\
	CC 2		& 0.4 &	2.5 	&		\\
	CC 3		& 0.1 &	0.5 	&		\\
	\midrule
	B 1		& &	4.0 &	Heating\\
	B 2		& &	3.0 &		\\
	B 3		& &	3.0 &		\\
	B 4 	& &	1.0	&		\\	
	\midrule
	EB 1	& 1.5 &	1.5 	&	Heating\\	
	EB 2	& 1.0 &	1.0		&		\\
	\midrule
	CHP 1	& 4.4 & 4	&	Heating, Electricity	\\
	CHP 2	& 2.2 &	2   &	 \\
	CHP 3	& 2.2 &	2   &	 \\
	CHP 4	& 1.1 &	1  	&		\\
	\bottomrule 
	\end{tabular}
\label{tab:components}
\end{table}

\clearpage

  \bibliographystyle{apalike}
  \renewcommand{\refname}{Bibliography}  
  \bibliography{literature.bib}

\end{document}